\documentclass[reqno,12pt,a4paper]{amsart}

\usepackage{mathrsfs}
\usepackage{amssymb}
\usepackage{amsfonts}
\usepackage{latexsym}
\usepackage{amsthm}

\usepackage{graphicx}
\def\lb{\label}

\newcommand{\er}[1]{\textrm{(\ref{#1})}}

\begin{document}


\renewcommand{\theequation}{\arabic{section}.\arabic{equation}}
\theoremstyle{plain}
\newtheorem{theorem}{\bf Theorem}[section]
\newtheorem{lemma}[theorem]{\bf Lemma}
\newtheorem{corollary}[theorem]{\bf Corollary}
\newtheorem{proposition}[theorem]{\bf Proposition}
\newtheorem{definition}[theorem]{\bf Definition}
\newtheorem{remark}[theorem]{\it Remark}

\def\a{\alpha}  \def\cA{{\mathcal A}}     \def\bA{{\bf A}}  \def\mA{{\mathscr A}}
\def\b{\beta}   \def\cB{{\mathcal B}}     \def\bB{{\bf B}}  \def\mB{{\mathscr B}}
\def\g{\gamma}  \def\cC{{\mathcal C}}     \def\bC{{\bf C}}  \def\mC{{\mathscr C}}
\def\G{\Gamma}  \def\cD{{\mathcal D}}     \def\bD{{\bf D}}  \def\mD{{\mathscr D}}
\def\d{\delta}  \def\cE{{\mathcal E}}     \def\bE{{\bf E}}  \def\mE{{\mathscr E}}
\def\D{\Delta}  \def\cF{{\mathcal F}}     \def\bF{{\bf F}}  \def\mF{{\mathscr F}}
\def\c{\chi}    \def\cG{{\mathcal G}}     \def\bG{{\bf G}}  \def\mG{{\mathscr G}}
\def\z{\zeta}   \def\cH{{\mathcal H}}     \def\bH{{\bf H}}  \def\mH{{\mathscr H}}
\def\e{\eta}    \def\cI{{\mathcal I}}     \def\bI{{\bf I}}  \def\mI{{\mathscr I}}
\def\p{\psi}    \def\cJ{{\mathcal J}}     \def\bJ{{\bf J}}  \def\mJ{{\mathscr J}}
\def\vT{\Theta} \def\cK{{\mathcal K}}     \def\bK{{\bf K}}  \def\mK{{\mathscr K}}
\def\k{\kappa}  \def\cL{{\mathcal L}}     \def\bL{{\bf L}}  \def\mL{{\mathscr L}}
\def\l{\lambda} \def\cM{{\mathcal M}}     \def\bM{{\bf M}}  \def\mM{{\mathscr M}}
\def\L{\Lambda} \def\cN{{\mathcal N}}     \def\bN{{\bf N}}  \def\mN{{\mathscr N}}
\def\m{\mu}     \def\cO{{\mathcal O}}     \def\bO{{\bf O}}  \def\mO{{\mathscr O}}
\def\n{\nu}     \def\cP{{\mathcal P}}     \def\bP{{\bf P}}  \def\mP{{\mathscr P}}
\def\r{\rho}    \def\cQ{{\mathcal Q}}     \def\bQ{{\bf Q}}  \def\mQ{{\mathscr Q}}
\def\s{\sigma}  \def\cR{{\mathcal R}}     \def\bR{{\bf R}}  \def\mR{{\mathscr R}}
\def\S{\Sigma}  \def\cS{{\mathcal S}}     \def\bS{{\bf S}}  \def\mS{{\mathscr S}}
\def\t{\tau}    \def\cT{{\mathcal T}}     \def\bT{{\bf T}}  \def\mT{{\mathscr T}}
\def\f{\phi}    \def\cU{{\mathcal U}}     \def\bU{{\bf U}}  \def\mU{{\mathscr U}}
\def\F{\Phi}    \def\cV{{\mathcal V}}     \def\bV{{\bf V}}  \def\mV{{\mathscr V}}
\def\P{\Psi}    \def\cW{{\mathcal W}}     \def\bW{{\bf W}}  \def\mW{{\mathscr W}}
\def\o{\omega}  \def\cX{{\mathcal X}}     \def\bX{{\bf X}}  \def\mX{{\mathscr X}}
\def\x{\xi}     \def\cY{{\mathcal Y}}     \def\bY{{\bf Y}}  \def\mY{{\mathscr Y}}
\def\X{\Xi}     \def\cZ{{\mathcal Z}}     \def\bZ{{\bf Z}}  \def\mZ{{\mathscr Z}}
\def\O{\Omega}

\newcommand{\gA}{\mathfrak{A}}
\newcommand{\gB}{\mathfrak{B}}
\newcommand{\gC}{\mathfrak{C}}
\newcommand{\gD}{\mathfrak{D}}
\newcommand{\gE}{\mathfrak{E}}
\newcommand{\gF}{\mathfrak{F}}
\newcommand{\gG}{\mathfrak{G}}
\newcommand{\gH}{\mathfrak{H}}
\newcommand{\gI}{\mathfrak{I}}
\newcommand{\gJ}{\mathfrak{J}}
\newcommand{\gK}{\mathfrak{K}}
\newcommand{\gL}{\mathfrak{L}}
\newcommand{\gM}{\mathfrak{M}}
\newcommand{\gN}{\mathfrak{N}}
\newcommand{\gO}{\mathfrak{O}}
\newcommand{\gP}{\mathfrak{P}}
\newcommand{\gR}{\mathfrak{R}}
\newcommand{\gS}{\mathfrak{S}}
\newcommand{\gT}{\mathfrak{T}}
\newcommand{\gU}{\mathfrak{U}}
\newcommand{\gV}{\mathfrak{V}}
\newcommand{\gW}{\mathfrak{W}}
\newcommand{\gX}{\mathfrak{X}}
\newcommand{\gY}{\mathfrak{Y}}
\newcommand{\gZ}{\mathfrak{Z}}

\def\ve{\varepsilon}   \def\vt{\vartheta}    \def\vp{\varphi}    \def\vk{\varkappa}

\def\Z{{\mathbb Z}}    \def\R{{\mathbb R}}   \def\C{{\mathbb C}}
\def\T{{\mathbb T}}    \def\N{{\mathbb N}}   \def\dD{{\mathbb D}}


\def\la{\leftarrow}              \def\ra{\rightarrow}            \def\Ra{\Rightarrow}
\def\ua{\uparrow}                \def\da{\downarrow}
\def\lra{\leftrightarrow}        \def\Lra{\Leftrightarrow}


\def\lt{\biggl}                  \def\rt{\biggr}
\def\ol{\overline}               \def\wt{\widetilde}
\def\no{\noindent}


\let\ge\geqslant                 \let\le\leqslant
\def\lan{\langle}                \def\ran{\rangle}
\def\/{\over}                    \def\iy{\infty}
\def\sm{\setminus}               \def\es{\emptyset}
\def\ss{\subset}                 \def\ts{\times}
\def\pa{\partial}                \def\os{\oplus}
\def\om{\ominus}                 \def\ev{\equiv}
\def\iint{\int\!\!\!\int}        \def\iintt{\mathop{\int\!\!\int\!\!\dots\!\!\int}\limits}
\def\el2{\ell^{\,2}}             \def\1{1\!\!1}
\def\sh{\sharp}
\def\wh{\widehat}

\def\where{\mathop{\mathrm{where}}\nolimits}
\def\as{\mathop{\mathrm{as}}\nolimits}
\def\Area{\mathop{\mathrm{Area}}\nolimits}
\def\arg{\mathop{\mathrm{arg}}\nolimits}
\def\const{\mathop{\mathrm{const}}\nolimits}
\def\det{\mathop{\mathrm{det}}\nolimits}
\def\diag{\mathop{\mathrm{diag}}\nolimits}
\def\diam{\mathop{\mathrm{diam}}\nolimits}
\def\dim{\mathop{\mathrm{dim}}\nolimits}
\def\dist{\mathop{\mathrm{dist}}\nolimits}
\def\Im{\mathop{\mathrm{Im}}\nolimits}
\def\Iso{\mathop{\mathrm{Iso}}\nolimits}
\def\Ker{\mathop{\mathrm{Ker}}\nolimits}
\def\Lip{\mathop{\mathrm{Lip}}\nolimits}
\def\rank{\mathop{\mathrm{rank}}\limits}
\def\Ran{\mathop{\mathrm{Ran}}\nolimits}
\def\Re{\mathop{\mathrm{Re}}\nolimits}
\def\Res{\mathop{\mathrm{Res}}\nolimits}
\def\res{\mathop{\mathrm{res}}\limits}
\def\sign{\mathop{\mathrm{sign}}\nolimits}
\def\span{\mathop{\mathrm{span}}\nolimits}
\def\supp{\mathop{\mathrm{supp}}\nolimits}
\def\Tr{\mathop{\mathrm{Tr}}\nolimits}
\def\BBox{\hspace{1mm}\vrule height6pt width5.5pt depth0pt \hspace{6pt}}


\newcommand\nh[2]{\widehat{#1}\vphantom{#1}^{(#2)}}
\def\dia{\diamond}

\def\Oplus{\bigoplus\nolimits}



\def\qqq{\qquad}
\def\qq{\quad}
\let\ge\geqslant
\let\le\leqslant
\let\geq\geqslant
\let\leq\leqslant
\newcommand{\ca}{\begin{cases}}
\newcommand{\ac}{\end{cases}}
\newcommand{\ma}{\begin{pmatrix}}
\newcommand{\am}{\end{pmatrix}}
\renewcommand{\[}{\begin{equation}}
\renewcommand{\]}{\end{equation}}
\def\bu{\bullet}

\title[{Resonance theory for perturbed Hill operator}]
{Resonance theory for perturbed Hill operator}

\date{\today}
\author[Evgeny Korotyaev]{Evgeny Korotyaev}
\address{Saint-Petersburg University, Russia
 \ korotyaev@gmail.com}


\begin{abstract}
We consider the Schr\"odinger operator $Hy=-y''+(p+q)y$ with a
periodic potential $p$  plus a compactly supported potential $q$ on
the real line. The spectrum of $H$ consists of an absolutely
continuous part plus a finite number of simple eigenvalues below the
spectrum and in each spectral gap $\g_n\ne \es, n\ge1$. We prove the
following results: 1)  the distribution of resonances in the disk
with large radius is determined, 2) the asymptotics of eigenvalues
and antibound states are determined at high energy gaps, 3) if $H$
has infinitely many open gaps in the continuous spectrum, then for
any sequence $(\vk)_1^\iy, \vk_n\in \{0,2\}$, there exists a
compactly supported potential $q$ with $\int_\R qdx=0$ such that $H$
has $\vk_n$ eigenvalues and $2-\vk_n$ antibound states (resonances)
in each gap $\g_n$ for $n$ large enough.

\end{abstract}

\maketitle

\section{Introduction}

Consider the Schr\"odinger operator $H$ acting in $L^2(\R)$ and given by
$$
H=H_0+q, \qqq {\rm where }\qqq H_0=-{d^2\/dx^2}+p.
$$
We assume that  $p\in L^2(0,1)$ is a real 1-periodic potential, and
$q$ is a real compactly supported potential and  belongs to
the class $\cQ_{t}^r$ given by
$$
\cQ_{t}^r=\rt\{q\in L^r(\R ): \  \supp q \ss [0,t]\rt\}, \qq t>0,\qq r\ge 1.
$$
The spectrum of $H_0$ is absolutely continuous and consists of
spectral bands $\gS_n$ separated by gaps $\g_n$, which are   given
by (see Fig. 1)
$$
\s(H_0)=\s_{ac}(H_0)=\cup_{n\ge 1} \gS_n,\qq
$$$$
\gS_n=[E^+_{n-1},E^-_n], \qq \g_{n}=(E^-_{n},E^+_n),\qq n\ge 1,\qq
{\rm and} \qq E_0^+<..\le E^+_{n-1}< E^-_n \le E^+_{n}<... \ \ .
$$
We assume that $E_0^+=0$.
The bands $\gS_n, \gS_{n+1}$ are separated by  a gap $\g_{n}=(E^-_{n},E^+_n)$.
If a gap degenerates, that is $\g_n=\es $, then the corresponding bands
$\gS_{n} $ and $\gS_{n+1}$ merge. Here $E_n^\pm$ is the eigenvalue of
the boundary value problem
\[
\lb{1}
-y''+p(x)y=\l y, \ \ \ \ \l\in \C ,\qqq y(x+2)=y(x), x\in \R.
\]
If $E_n^-=E_n^+$ for some $n$, then this number $E_n^{\pm}$ is the
double eigenvalue of the problem \er{1}. The lowest eigenvalue
$E_0^+=0$ is always simple  and the corresponding eigenfunction is
1-periodic. The eigenfunctions, corresponding to the eigenvalue
$E_{2n}^{\pm}$, are 1-periodic, and for the case $E_{2n+1}^{\pm}$
they are anti-periodic,  i.e., $y(x+1)=-y(x),\ \ x\in\R$.

It is well known, that the spectrum of $H$ consists of
an absolutely continuous part $\s_{ac}(H)=\s(H_0)$ plus
a finite number of simple eigenvalues, both in each gap $\g_n\ne \es, n\geq 1$
 and in the half-line $(-\iy,E_0^+)$, see  \cite{Rb}, \cite{F1}.
 Moreover, in a remote open gap $\g_n$  the operator $H$ has
 at most two  eigenvalues \cite{Rb} and precisely one  eigenvalue
 in the case $\int_\R q(x)dx\ne 0$ \cite{Zh1}, \cite{F2}.
Note that  the potential $q$ in \cite{Rb},\cite{Zh1} belongs to the
more general class, see also \cite{F3}-\cite{F4}, \cite{GS},
\cite{So}, \cite{Zh2}, \cite{Zh3}.

The resonance theory for the multidimensional Schr\"odinger operator
 with a periodic
potential plus a real compactly supported potential
 has been
much less studied, see \cite{D}, \cite{G} and references therein.
Some results for the case of a slowly varying perturbations of a 1D
periodic Schr\"odinger operator have been announced in \cite{KM}.

\begin{figure}
\tiny
\unitlength=1.00mm
\special{em:linewidth 0.4pt}
\linethickness{0.4pt}
\begin{picture}(108.67,33.67)
\put(41.00,17.33){\line(1,0){67.67}}
\put(44.33,9.00){\line(0,1){24.67}}
\put(108.33,14.00){\makebox(0,0)[cc]{$\Re\l$}}
\put(41.66,33.67){\makebox(0,0)[cc]{$\Im\l$}}
\put(42.00,14.33){\makebox(0,0)[cc]{$0$}}
\put(44.33,17.33){\linethickness{4.0pt}\line(1,0){11.33}}
\put(66.66,17.33){\linethickness{4.0pt}\line(1,0){11.67}}
\put(82.00,17.33){\linethickness{4.0pt}\line(1,0){12.00}}
\put(95.66,17.33){\linethickness{4.0pt}\line(1,0){11.00}}
\put(46.66,20.00){\makebox(0,0)[cc]{$E_0^+$}}
\put(56.66,20.33){\makebox(0,0)[cc]{$E_1^-$}}
\put(68.66,20.33){\makebox(0,0)[cc]{$E_1^+$}}
\put(78.33,20.33){\makebox(0,0)[cc]{$E_2^-$}}
\put(84.33,20.33){\makebox(0,0)[cc]{$E_2^+$}}
\put(93.00,20.33){\makebox(0,0)[cc]{$E_3^-$}}
\put(98.66,20.33){\makebox(0,0)[cc]{$E_3^+$}}
\put(106.33,20.33){\makebox(0,0)[cc]{$E_4^-$}}
\end{picture}
\caption{The cut domain $\C\sm \cup \gS_n$ and the cuts (bands)
$\gS_n=[E^+_{n-1},E^-_n], n\ge 1$}
\lb{sS}
\end{figure}
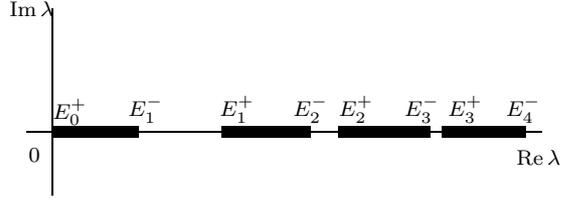

Introduce the two-sheeted Riemann surface $\L$  obtained by joining
the upper and lower rims of two copies of the cut plane
$\C\sm\s_{ac}(H_0)$ in the usual (crosswise) way. The n-th gap on
the first physical sheet $\L_1$ we will denote by $\g_n^{(1)}$ and
whereas  the same gap on the second nonphysical sheet $\L_2$ we will
denote by $\g_n^{(2)}$. Let $\g_n^c$ be the union of $\ol\g_n^{(1)}$
and $\ol\g_n^{(2)}$, i.e.,
 $$
 \g_n^c=\ol\g_n^{(1)}\cup \ol\g_n^{(2)}.
 $$
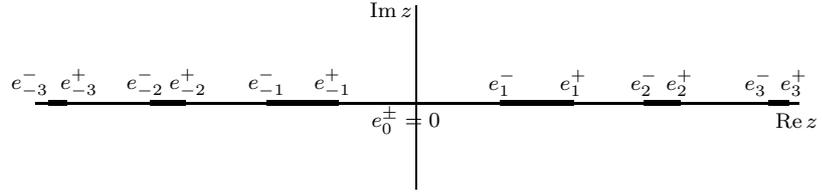
\begin{figure}
\tiny
\unitlength=1mm
\special{em:linewidth 0.4pt}
\linethickness{0.4pt}
\begin{picture}(120.67,34.33)
\put(20.33,21.33){\line(1,0){100.33}}
\put(70.33,10.00){\line(0,1){24.33}}
\put(69.00,19.00){\makebox(0,0)[cc]{$e_0^\pm=0$}}
\put(120.33,19.00){\makebox(0,0)[cc]{$\Re z$}}
\put(67.00,33.67){\makebox(0,0)[cc]{$\Im z$}}
\put(81.33,21.33){\linethickness{2.0pt}\line(1,0){9.67}}
\put(100.33,21.33){\linethickness{2.0pt}\line(1,0){4.67}}
\put(116.67,21.33){\linethickness{2.0pt}\line(1,0){2.67}}
\put(60.00,21.33){\linethickness{2.0pt}\line(-1,0){9.33}}
\put(40.00,21.33){\linethickness{2.0pt}\line(-1,0){4.67}}
\put(24.33,21.33){\linethickness{2.0pt}\line(-1,0){2.33}}
\put(81.67,24.00){\makebox(0,0)[cc]{$e_1^-$}}
\put(91.00,24.00){\makebox(0,0)[cc]{$e_1^+$}}
\put(100.33,24.00){\makebox(0,0)[cc]{$e_2^-$}}
\put(105.00,24.00){\makebox(0,0)[cc]{$e_2^+$}}
\put(115.33,24.00){\makebox(0,0)[cc]{$e_3^-$}}
\put(120.00,24.00){\makebox(0,0)[cc]{$e_3^+$}}
\put(59.33,24.00){\makebox(0,0)[cc]{$e_{-1}^+$}}
\put(50.67,24.00){\makebox(0,0)[cc]{$e_{-1}^-$}}
\put(40.33,24.00){\makebox(0,0)[cc]{$e_{-2}^+$}}
\put(34.67,24.00){\makebox(0,0)[cc]{$e_{-2}^-$}}
\put(26.00,24.00){\makebox(0,0)[cc]{$e_{-3}^+$}}
\put(19.50,24.00){\makebox(0,0)[cc]{$e_{-3}^-$}}
\end{picture}
\caption{The cut domain $\cZ=\C\sm \cup \ol g_n$ and
the cuts $g_n=(e_n^-,e_n^+)$ in the $z$-plane.}
\lb{z}
\end{figure}
In what follows we will use the momentum variable $z=\sqrt \l, \l\in
\L$ and the corresponding Riemann surface $\cM$, which is more
convenient for us, than the Riemann surface $\L$. The mapping
$\l\mapsto z=\sqrt \l$ is a bijection between the  cut Riemann
surface $\L\sm \cup \g_n^c$ and the cut  momentum  domain $\cZ$
(see Fig.2.)  given by
 \[
\lb{Z0}
\cZ=\C\sm \bigcup_{n\ne 0} \ol g_n, \qq {\rm where} \ \
g_n=(e_n^-,e_n^+), \qq e_n^\pm=-e_{-n}^\mp= \sqrt{E_n^\pm}>0,\qq
n\ge 1.
\]
Here $\R\sm  \bigcup_{n\ne 0} \ol g_n$ is the momentum spectrum and
$g_n\ne \es$ is the momentum gap. Slitting the n-th momentum gap
$g_n$ (suppose it is nontrivial) we obtain a cut $g_n^c$ with an
upper rim  $g_n^+$ and lower rim  $g_n^-$. Below we will identify
this cut $g_n^c$ and the union of of the upper rim  (gap) $\ol
g_{n}^+$ and the lower rim (gap) $\ol g_{n}^{\ -}$, i.e.,
\[
g_n^c=\ol g_{n}^+\cup \ol g_{n}^-.
\]
In order to construct the  Riemann surface $\cM$ we take the cut
domain $\cZ$ and identify (i.e. we glue) the upper rim $g_{n}^+$ of
the  cut $g_n^c$ with the upper rim $g_{-n}^+$  of the cut
$g_{-n}^c$ and correspondingly the lower rim $g_{n}^-$ of the  cut
$g_n^c$ with the  lower  rim $g_{-n}^-$  of the  cut $g_{-n}^c$ for
all nontrivial gaps. The mapping $\l\mapsto z=\sqrt \l$ from $\L$
onto $\cM$ is one-to-one and onto and we have the following.


1) The physical  gap $\g_n^{(1)}\ss \L_1$ is mapped onto
the physical "gap" (the upper rim) $g_n^+\ss \cM_1$
and the half-line $(-\iy,0)\ss \L_1$ is mapped onto $i\R_+$.

2) The nonphysical gap $\g_n^{(2)}\ss \L_2$ is mapped onto
the nonphysical "gap" (the lower rim) $g_n^-\ss \cM_2$ and
the half-line  $(-\iy, 0)\ss \L_2$  is mapped onto  $i\R_-$.

3) $\C _+=\{z:\Im z>0\}$ plus all physical gaps $g_{n}^+$ is
a so-called physical "sheet" $\cM_1$.

4) $\C _-=\{z:\Im z<0\}$ plus all nonphysical gaps $g_{n}^-$
is a so-called nonphysical "sheet" $\cM_2$.

5) The momentum spectrum $\s_M=\R\sm \cup [e_n^-,e_n^+]$ joints the
first  $\cM_1$ and second sheets  $\cM_2$.
\medskip

Note that if $p=0$, then  $\L$ is a Riemann surface  of the function
$\sqrt \l$,  $\cM=\C$ is  the momentum plane,  $\cM_1=\C_+$ is the
physical "sheet" and  $\cM_2=\C_-$ is the nonphysical "sheet".

We introduce the determinant
$$
D(z)=\det (I+q(H_0-z^2)^{-1}),\qqq z\in \C_+,\qqq
$$
which is analytic in $\C_+$  and continuous up to $\R\sm \{z:
z=e_n^\pm, n\in \Z \}$, where $e_0^\pm=0$,  see \cite{F4},
\cite{F1}. It is well known that if $D(z)=0$ for some zero $z\in
\cM_1$,  then $z^2$ is an eigenvalue of $H$ and $z\in\cup_{n\ne0}
g_n^+$ or $z\in i\R_+$.  We introduce our basic function $\x$ by
\[
\x(z)=2i\sin k(z)D(z), \qqq z\in\C_+.
\]
Here $k(z)$ is the quasimomentum for the operator $H_0$ introduced
by Firsova \cite{F3} and Marchenko-Ostrovski \cite{MO}, see Section
2 for a precise definition of $k$. In Section 2 we describe the
properties of the function $k$, which is analytic  in $z\in\cZ$.
Moreover, we show that $\sin k(z), z\in \cZ$  is analytic in $\cM$
and $\cM$ is the Riemann surface of $\sin k(z)$. All zeros of $\sin
k(z), z\in \cM$ have the form $e_n^\pm, n\in \Z$, where $e_0^\pm=0$.
In Theorem \ref{T1} we will show that $\x$ has an analytic extension
from $\C_+$ into the Riemann surface $\cM$.

\no {\bf Definition S.}
Let $\z\in \cM$ be a zero of $\x(z), z\in \cM$ and assume that
$\z\ne e_n^+$ for any $e_n^+=e_n^-, n\ne 0$.

\no 1) If  $\z\in i\R_+$ or $\z\in\bigcup\limits_{n\ne 0} g_n^+$,
we call $\z$ a {\it bound state}.

\no 2) If  $\z\in \cM_2$ and $\z\ne e_n^\pm, n\in \Z$,
we call $\z$ a {\it resonance}.

\no 3) Let  $e_0^\pm=0$.  If  $\z=0$ or $\z=e_n^\pm$
for the open gap $|g_n|>0, n\ne 0$, we call $\z$ a {\it  virtual state}.

\no 4) A point $\z\in\cM$ is called a {\it state} if it is either a bound
state or a resonance or a virtual state.
We denote by $\gS_{st}(H)$ the set of all states.
If $\z\in g_n^-, n\ne 0$ or $\z\in i\R_-$, then
we call $\z$ an {\it antibound state}.

\no 5) The {\it multiplicity} of a bound state, a resonance or the
point $0$ is the multiplicity of the corresponding zero. The {\it
multiplicity} of the virtual state $\z\ne 0$ is the multiplicity of the
zero $z=0$ of the function $\x(\z+z^2)$. A state with multiplicity
one is called simple.

\bigskip

Of course, $z^2$ is really the energy, but since
the momentum $z$ is the natural parameter, we will abuse the terminology.

We recall the results about the resonances from \cite{F1}:

1) Let $q=0$. Thus we have $D=1$ and $\x(z)=2i\sin k(z), \ z\in\cZ$.
Then the operator $H_0$ has only virtual states $e_n^\pm$
for all $e_n^-\ne e_n^+, n\ne 0$ and
$e_0^+=0$. There are no other states.

2) Let a gap  $g_n=\es$ for  some $n\ne 0$. Then for any  $h\in
C_0^\iy(\R), h\ne 0$ the  function $((H-z^2)^{-1}h,h)$ is analytic
at the point $e_n^+=e_n^-\in \cZ$. Roughly speaking there is no
difference between these points and other points inside the spectrum
$\s_{ac}(H)$. The point $e_n^+=e_n^-$ is not the state.

3) If $\int_\R q(x)dx\ne 0$, then $H$ has  precisely one bound state
on each open physical gap and an odd number $\ge 1$ of antibound
states on each open non-physical gap for $n$ large enough.

Define the coefficients for all $n\ge 1$ by
\[
\lb{fco}
\wh q_0=\int_\R q(x)dx,\ \
\wh q_n=\wh q_{cn}+i\wh q_{sn},\ \
\wh q_{cn}=\int_\R q(x)\cos 2\pi nxdx,\ \
\wh q_{sn}=\int_\R q(x)\sin 2\pi nxdx.
\]

Let $\m_n^2,  n\ge 1$ be eigenvalues  and $y_n$ be
the corresponding eigenfunctions of the Sturm-Liouville problem
\[
\lb{Dei}
-y_n''+py_n=\m_n^2y_n ,\qqq \qqq y_n(0)=y_n(1)=0,
\]
on the interval $[0,1]$. It is well known that each $\m_n^2\in
[E^-_n,E^+_n ]$ for all $n\ge 1$.

In order to formulate Theorem \ref{T1} we define $c_n, s_n$ the angles $\f_n\in [0,2\pi )$
\[
\lb{csf}
c_n=\cos \f_n, \qq
s_n=\sin \f_n\in [-1,1],\qqq \f_n\in [0,2\pi )
\]
by the identities
\[
\lb{cs} {E_n^-+E_n^+\/2}-\m_n^2={|\g_n|\/2}c_n, \qq
|1-c_n^2|^{1\/2}\sign (|y_n'(1)|-1)=s_n, \ \ \ {\rm if} \qq
|\g_n|>0,
\]
where  all  eigenfunctions satisfy $y_n'(0)=1$. We describe our
first main results about states.

\begin{theorem}
\lb{T1}
Let potentials $(p,q)\in L^2(0,1)\ts \cQ_t^1, t>0$. Then we have

i) $\x$ has an analytic extension from $\C_+$ into the Riemann
surface $\cM$ and the function $J(z)=\Re \x(z), z\in \s_M=\R\sm \cup
[e_n^-,e_n^+]$ has an analytic extension into the whole  plane $\C$.

ii) There exist an even number $\ge 0$ of states (counted with multiplicity)
  on each set $g_n^c\ne \es,n\ne 0$, where $g_n^c$ is a union of
  the physical  gap $\ol g_n ^{\ +}\ss \cM_1$ and
  the non-physical gap $\ol g_n^{\ -}\ss\cM_2$.

 iii) There are no states in the "forbidden"  domain $\cD\ss\C_-$ given by
\begin{multline}
\lb{T1-3}
\cD=\rt\{z\in \C_-: |z|>\max \{180e^{2\|p\|_1},\ C_0e^{2t|\Im z|}\}\rt\},\qq
C_0=12\|q\|_te^{\|p\|_1+\|q\|_t+2\|p\|_t},
\end{multline}
where $\|q\|_t=\!\int_0^t\!\!|q(x)|dx$.

 iv)  In each  $g_n^c\ne \es, n\ge 1+{e^{t\pi /2}\/\pi}C_0$
 there exist exactly two   simple real states $z_n^\pm\in g_n^c$
 such that $e_n^-\le z_n^-<e_n< z_n^+\le e_n^+$
 (for the definition of $e_n$ see \er{dzn}) and satisfy
\begin{multline}
\lb{T2-1ad}
z_n^\pm=e_n^\pm\mp {2|\g_n|\/(4\pi n)^3}
(\mp \wh q_0-c_n\wh q_{cn}+s_n\wh q_{sn}+ O(1/n))^2,\qqq
\\
\qqq (-1)^{n}J(z_n^\pm)={|\g_n|\/(2\pi n)^2}(\mp \wh q_0-c_n\wh
q_{cn}+s_n\wh q_{sn}+ O(1/n))\qqq\as \qqq n\to \iy.
\end{multline}
Moreover, if a state $\z\in \{z_n^-, z_n^+\}$ satisfies
$(-1)^{n}J(\z)>0$, $(-1)^{n}J(\z)<0$, $J(\z)=0$), then $\z$ is a
bound state,  an antibound state or a virtual state correspondingly
and, in particular,
\begin{multline}
\lb{q0} if \qq \wh q_0>0 \ \Rightarrow \qq z_n^-\in \cM_1\ is \
bound  \ state, \qqq
z_n^+\in \cM_2\ is \ antibound  \ state, \\
 if \qq \wh q_0<0 \ \Rightarrow \ \qq
z_n^-\in \cM_2 \ is \ antibound  \ state, \qqq  z_n^+\in \cM_1\ is \
bound  \ state.
\end{multline}
\end{theorem}

\no {\bf Remark.} 1) The forbidden domain $\cD$  is similar to
the case $p=0$, see \cite{K2} and Fig. 3.

2) In the proof of Theorem
\ref{T2} the estimates $z_n^-<e_n< z_n^+$ and asymptotics \er{T2-1ad}
are important.
\begin{figure}
\tiny
\unitlength 1mm 
\linethickness{0.4pt}
\ifx\plotpoint\undefined\newsavebox{\plotpoint}\fi 
\begin{picture}(76.222,83.464)(0,0)
\put(-1.842,54.344){\line(1,0){78.064}}
\put(1.854,72.6){\line(0,-1){69.664}}
\qbezier(67.71,54.344)(69.782,56.136)(72.078,54.344)
\qbezier(53.71,54.344)(55.782,56.136)(58.078,54.344)
\qbezier(67.822,54.344)(69.39,53.056)(71.854,54.232)
\qbezier(53.822,54.344)(55.39,53.056)(57.854,54.232)
\qbezier(11.486,54.232)(16.974,56.696)(21.118,54.232)
\qbezier(11.486,54.232)(16.19,52.72)(20.894,54.12)
\qbezier(31.31,54.456)(34.95,56.696)(39.71,54.456)
\qbezier(30.974,54.232)(35.006,52.44)(39.486,54.232)
\bezier{50}(1.854,9.88)(21.622,10.888)(34.446,24.44)
\bezier{50}(34.446,24.328)(45.534,36.816)(45.87,54.456)
\qbezier(21.342,48.968)(35.398,24.776)(75.662,19.848)
\multiput(16.19,38.888)(.0331852,.0373333){27}{\line(0,1){.0373333}}
\multiput(19.662,43.144)(.0331852,.0373333){27}{\line(0,1){.0373333}}
\multiput(34.222,55.016)(.0331852,.0373333){27}{\line(0,1){.0373333}}
\multiput(54.606,54.344)(.0331852,.0373333){27}{\line(0,1){.0373333}}
\multiput(56.51,53.112)(.0331852,.0373333){27}{\line(0,1){.0373333}}
\multiput(70.622,54.568)(.0331852,.0373333){27}{\line(0,1){.0373333}}
\multiput(68.27,53.448)(.0331852,.0373333){27}{\line(0,1){.0373333}}
\multiput(37.134,54.904)(.0331852,.0373333){27}{\line(0,1){.0373333}}
\multiput(37.134,53.112)(.0331852,.0373333){27}{\line(0,1){.0373333}}
\multiput(32.766,52.888)(.0331852,.0373333){27}{\line(0,1){.0373333}}
\multiput(26.046,46.056)(.0331852,.0373333){27}{\line(0,1){.0373333}}
\multiput(1.406,59.944)(.0331852,.0373333){27}{\line(0,1){.0373333}}
\multiput(1.406,65.432)(.0331852,.0373333){27}{\line(0,1){.0373333}}
\multiput(1.406,69.576)(.0331852,.0373333){27}{\line(0,1){.0373333}}
\multiput(30.078,32.392)(.0331852,.0373333){27}{\line(0,1){.0373333}}
\multiput(44.19,24.216)(.0331852,.0373333){27}{\line(0,1){.0373333}}
\multiput(48.222,17.496)(.0331852,.0373333){27}{\line(0,1){.0373333}}
\multiput(57.518,21.416)(.0331852,.0373333){27}{\line(0,1){.0373333}}
\multiput(63.342,15.592)(.0331852,.0373333){27}{\line(0,1){.0373333}}
\multiput(71.07,18.952)(.0331852,.0373333){27}{\line(0,1){.0373333}}
\multiput(16.078,39.784)(.0373333,-.0331852){27}{\line(1,0){.0373333}}
\multiput(19.55,44.04)(.0373333,-.0331852){27}{\line(1,0){.0373333}}
\multiput(34.11,55.912)(.0373333,-.0331852){27}{\line(1,0){.0373333}}
\multiput(54.494,55.24)(.0373333,-.0331852){27}{\line(1,0){.0373333}}
\multiput(56.398,54.008)(.0373333,-.0331852){27}{\line(1,0){.0373333}}
\multiput(70.51,55.464)(.0373333,-.0331852){27}{\line(1,0){.0373333}}
\multiput(68.158,54.344)(.0373333,-.0331852){27}{\line(1,0){.0373333}}
\multiput(37.022,55.8)(.0373333,-.0331852){27}{\line(1,0){.0373333}}
\multiput(37.022,54.008)(.0373333,-.0331852){27}{\line(1,0){.0373333}}
\multiput(32.654,53.784)(.0373333,-.0331852){27}{\line(1,0){.0373333}}
\multiput(25.934,46.952)(.0373333,-.0331852){27}{\line(1,0){.0373333}}
\multiput(1.294,60.84)(.0373333,-.0331852){27}{\line(1,0){.0373333}}
\multiput(1.294,66.328)(.0373333,-.0331852){27}{\line(1,0){.0373333}}
\multiput(1.294,70.472)(.0373333,-.0331852){27}{\line(1,0){.0373333}}
\multiput(29.966,33.288)(.0373333,-.0331852){27}{\line(1,0){.0373333}}
\multiput(44.078,25.112)(.0373333,-.0331852){27}{\line(1,0){.0373333}}
\multiput(48.11,18.392)(.0373333,-.0331852){27}{\line(1,0){.0373333}}
\multiput(57.406,22.312)(.0373333,-.0331852){27}{\line(1,0){.0373333}}
\multiput(63.23,16.488)(.0373333,-.0331852){27}{\line(1,0){.0373333}}
\multiput(70.958,19.848)(.0373333,-.0331852){27}{\line(1,0){.0373333}}
\put(68.158,40.568){\makebox(0,0)[cc]{$the \ forbiden\ domain \
\cD$}} \put(0.158,74.568){\makebox(0,0)[cc]{$ \Im z$}}
\put(81.158,54.568){\makebox(0,0)[cc]{$ \Re z$}}
\put(67.262,72.488){\makebox(0,0)[cc]{$the \ domain \ \cZ \ $}}
\put(3.15,51.72){\makebox(0,0)[cc]{$0$}}
\put(16.158,58.568){\makebox(0,0)[cc]{$ g_1^+$}}
\put(16.158,51.568){\makebox(0,0)[cc]{$ g_1^-$}}
\put(36.158,58.568){\makebox(0,0)[cc]{$ g_2^+$}}
\put(36.158,51.568){\makebox(0,0)[cc]{$ g_2^-$}}
\put(55.158,58.568){\makebox(0,0)[cc]{$ g_3^+$}}
\put(55.158,51.568){\makebox(0,0)[cc]{$ g_3^-$}}
\put(69.158,58.568){\makebox(0,0)[cc]{$ g_4^+$}}
\put(69.158,51.568){\makebox(0,0)[cc]{$ g_4^-$}}
\end{picture}
\caption{\footnotesize The bound states and resonances on the $\cZ$
domain with the physical rims $g_n^+$ and the nonphysical  rims
$g_n^-$} \lb{res}
\end{figure}
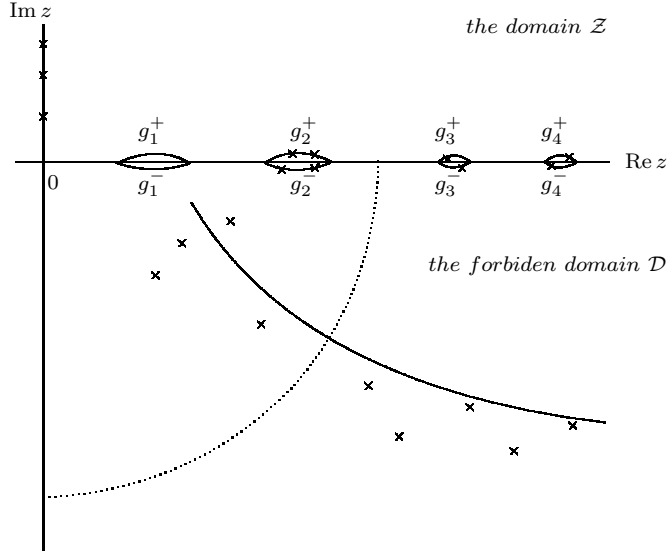

Recall that $p$ is even, i.e., $p\in L_{even}^2(0,1)=\{p\in
L^2(0,1), p(x)=p(1-x), x\in (0,1)\}$ iff $\m_n^2\in \{E_n^-,E_n^+\}$
for all $n\ge 1$, see \cite{GT}, \cite{KK1}. Note that $\m_n^2\in
\{E_n^-,E_n^+\} \Leftrightarrow s_n=0$.

\begin{corollary}
\lb{Tc} Let $\wh q_0=0$ and $\wh q_n=|\wh q_n|e^{i\t_n}, n\ge 1$ for
some $\t_n\in [0,2\pi)$. Assume that $|\cos (\f_n+\t_n)|>\ve>0$ and
$|\wh q_n|>n^{-\a}$ for  $n$ large enough and for some  $\ve, \a\in
(0,1)$, where $\f_n$ is defined by \er{cs}. Then we have

i) The operator $H$ has $\vk_n=1-\sign \cos (\f_n+\t_n)$ bound
states in the physical gap $g_n^+\ne \es$ and $2-\vk_n$ antibound
states inside the  nonphysical gap $g_n^-\ne \es$  for $n$ large
enough.

ii) Let in addition, a real potential $V\in \cQ_t^1$ and let $|\wh
V_{n}|=o(n^{-\a})$ as $n\to \iy$. Then  the operator $H+V$ has
$\vk_n$ bound states in the gap $g_n^+\ne \es$ and $2-\vk_n$
antibound states inside the gap $g_n^-\ne \es$  for $n$ large
enough.

iii) Let in addition, $p\in L_{even}^2(0,1)$.
Then the following asymptotics hold true:
\[
\lb{T4-1} z_n^\pm=e_n^\pm\mp {2|\g_n|\/(4\pi n)^3}\rt(c_n\wh q_{cn}+
{O(1)\/n}\rt)^2 \qqq \as \qq n\to \iy.
\]
Moreover,  if $|\wh q_{cn}|>n^{-\a}$ for   $n$ large enough, then
$H$ has exactly $\vk_n=1-\sign c_n\wh q_{cn}$ bound states in
each open gap $g_n^+$
and  $2-\vk_n$ antibound inside gap $g_n^-\ne \es$ for $n$ large enough.

\end{corollary}

{\bf Remark.} Let  all conditions in Corollary \ref{Tc} iii) hold
true and let $\wh q_{cn}>n^{-\a}$ for all $n$ large enough. Then

 if $\m_n=e_n^-$, then $H$ has exactly $2$ bound states in
 the open gap $g_n^+$ for $n$ large enough,

 if $\m_n=e_n^+$, then $H$ has not bound states in
 the open gap $g_n^+$ for $n$ large enough.

\bigskip

Let $\#(H,r, X)$ be the total number of state of $H$ in the set
$X\subseteq \cM$ having modulus $\le r$, each state being counted
according to its multiplicity.

\begin{theorem}
\lb{T2}
Let  the real potential $q\in L^2(\R )$ and let
$[0,t]$ be the convex hull of\ the support of $q$ for some $t>0$. Then
 the following asymptotics hold true:
\[
\lb{T2-2}
\#(H,r,\C_-)=r{2t+o(1)\/\pi}\qq as \qq r\to \iy.
\]
\end{theorem}

\no{\bf Remark}. 1) The first term in \er{T2-2} does not depend
on the periodic potential $p$.

\medskip

\no 2) The distribution of the resonances for the case $p\ne \const$
in the domain  $\C_-$
is similar to the case $p=0$, obtained by Zworski \cite{Z1}.

\medskip

\no 3) In the proof of \er{T2-2} we use the Paley-Wiener type Theorem
from \cite{Fr}, the Levinson Theorem
(see Sect. 5) and a priori estimates from \cite{KK}, \cite{MO}.

\medskip

Consider  some inverse problems for the operator $H$.

\begin{theorem}
\lb{T3} i)  Let the spectrum of the operator $H_0$ have infinitely
many gaps $\g_n\ne \es$ for some $p\in L^2(0,1)$. Then for any
sequence $(\vk_n)_{1}^\iy$, where $\vk_n\in \{0,2\}$,  there exists
some potential $q\in \cQ_1^1$ (defined by \er{defq}) such that $H$
has exactly $\vk_n$ bound states in each gap $g_n^+\ne \es$  and
$2-\vk_n$ antibound states inside each gap $g_n^-\ne \es$ for $n$
large enough.

\no ii)  Let  $q\in \cQ_t^1$ satisfy   $\wh q_0=0$ and let $|\wh
q_{n}|>n^{-\a}$ for all $n$ large enough and some $\a\in (0,1)$.
Then for any sequences $(\vk_n)_{1}^\iy$, where $\vk_n\in \{0,2\}$
and $(\d_n)_1^\iy\in \ell^2$, where all $\d_n\ge 0$ and infinitely
many $\d_n>0$, there exists a potential $p\in L^2(0,1)$  such that
each gap length  $|\g_n|=\d_n, n\ge 1$. Moreover, $H$ has exactly
$\vk_n$ bound states in each physical gap $g_n^+\ne \es$ and
$2-\vk_n$  antibound states  inside each non-physical gap $g_n^-\ne
\es$ for $n$ large enough.
\end{theorem}

{\bf Remark .}  The proof of ii) is more difficult and here we use
results from the inverse spectral theory from \cite{K5}.

A lot of papers are devoted to the resonances for the Schr\"odinger
operator with $p=0$, see \cite{Fr}, \cite{H}, \cite{K1}, \cite{K2},
\cite{S}, \cite{Z1}, \cite{Z3}  and references therein.
 Although resonances have been studied in many
settings, but there are relatively few cases where the asymptotics
of the resonance counting function are known, mainly the one
dimensional case \cite{Fr}, \cite{K1}, \cite{K2}, \cite{S}, and
\cite{Z1}. We recall that Zworski [Z1] obtained the first results
about the distribution of resonances for the Schr\"odinger operator
with compactly supported potentials on the real line. The author
obtained the characterization (plus uniqueness and recovering) of
$S$-matrix for the Schr\"odinger operator with a compactly supported
potential on the real line \cite{K2} and the half-line \cite{K1},
see also \cite{Z2}, \cite{BKW} about uniqueness.

For the Schr\"odinger operator on the half line the author \cite{K3}
obtained the following stability results.

(i) Let $z_n, n\ge 1$ be the sequence of all zeros (all states) of
the Jost function for some real compactly supported potential $q$.
Assume that $\sum_{n\ge 1}n^3|z_n-\wt z_n|^2<\iy$ for some sequence
$\wt z_n\in \C, n\ge 1$. Then $\wt z_n$ is the sequence of all
zeroes of the Jost function for some unique real compactly supported
potential.

(ii) The measure associated with the zeros of the Jost function is
the Carleson measure, and the sum $\sum (1+|z_n|)^{-\a}, \a>1$ is
estimated in terms of the $L^1$-norm of the potential $q$.

Brown and Weikard \cite{BW} considered  the Schr\"odinger operator
$-y''+(p_A+q)y$ on the half-line, where $p_A$ is an
algebro-geometric potentials and $q$ is a compactly supported
potential. They proved that the zeros of the Jost function determine
$q$ uniquely.

Christiansen \cite{Ch} considered resonances associated to the Schr\"odinger
operator $-y''+(p_{S}+q)y$ on the real line, where $p_S$ is a step potential. She determined asymptotics of the resonance-counting function. Moreover, she proved that the
resonances determine $q$ uniquely.

The plan of the paper is as follows. In Section 2 we describe the
preliminary results about fundamental solutions for the operator
$H_0$. In Section 3 we describe the scattering for $H, H_0$. In
Sections 4 we study the function $\x$ and the entire function
$F=\z(z)\z(-z)$. Here it is important that $\x$ has a finite number
zeros in $\C_+$. That makes possible to reformulate the problem for
the differential operator theory as a problem of the entire function
theory and  the conformal mapping theory. Thus we should study the
function $F$ using various "geometric properties" of conformal
mappings from \cite{KK}, \cite{MO}. The properties of $F$ are the
key of the proof of main Theorems \ref{T1}-\ref{T3}, given in
Section 5.

\section {The unperturbed operator $H_0$ }
\setcounter{equation}{0}

{\bf 2.1. Fundamental solutions.} In order to describe the spectral
properties of the operator $H_0$, we start from the properties of
the canonical fundamental system $\vt, \vp$ of the equation
$-y''+py=z^2y, z\in \C$ with the initial conditions
$\vp'(0,z)=\vt(0,z)=1$ and $\vp(0,z)=\vt'(0,z)=0$, where $u'=\pa_x
u$. They satisfy the integral equations
\begin{multline}
\lb{fs}
\vt(x,z)=\cos zx+\int_0^x{\sin z(x-s)\/z}p(s)\vt(s,z)ds,\qqq  \\
\vp(x,z)={\sin zx\/z}+\int_0^x{\sin z(x-s)\/z}p(s)\vp(s,z)ds.
\end{multline}
For each $x\in \R$ the functions $\vt (x,z), \vp (x,z)$
are entire in $z\in\C$ and satisfy
\begin{multline}
\lb{efs1}
\max \{|z|_1|\vp(x,z)|, \ |\vp'(x,z)| , |\vt(x,z)|,
{1\/|z|_1}|\vt'(x,z)|    \} \le X,\\
|\vp(x,z)-{\sin zx\/z}|\le {X\/|z|_1^2}\|p\|_x,
\qq |\vt(x,z)-{\cos zx}|\le {X\/|z|_1}\|p\|_x,
\end{multline}
for all $(p,x,z)\in L_{loc}^1(\R)\ts \R\ts \C$, see p. 13 in  \cite{PT}, where
$$
 X=e^{|\Im z|x+\|p\|_x}, \qqq \|p\|_t=\int_0^t|p(s)|ds,\qqq |z|_1=\max \{1,|z|\}.
$$

Let $\vp(x,z,\t)$   be the solutions of the equation with
the parameter $\t\in \R$
\[
\lb{x+t}
-\vp''+p(x+\t)\vp=z^2 \vp, \qq \ \vp(0,z,\t)=0,\qq \vp'(0,z,\t)=1, \qq  z\in \C.
\]
This solution  $\vp(x,z,\t)$ is expressed in terms of $\vt(x,\cdot),
\vp(x,\cdot)$ by
\[
\lb{ffs} \vp(x,\cdot,\t)
=\vt(\t,\cdot)\vp(x+\t,\cdot)-\vp(\t,\cdot)\vt(x+\t,\cdot).
\]
The function $\vp(1,z,x)$ for all $(x,z)\in \R\ts \C$
satisfies the following identity (see \cite{Tr})
\[
\lb{if}
\vp(1,\cdot,\cdot)=\vp(1,\cdot)\vt^2-\vt'(1,\cdot)\vp^2+2\b\vp\vt=
\vp(1,\cdot)\p_-\p_+.
\]

Below we need a solution of the equation
$-y''+(p-z^2)y=f, \ y(0)=y'(0)=0$  given by
\[
\lb{ieq}
y=\int_0^x\vp(x-\t,z,\t)f(\t)d\t.
\]
Recall that $\m_n^2$ is the Dirichlet eigenvalue, defined by \er{Dei}.
 The following asymptotics hold true as $n\to \iy$ (see \cite{PT}, \cite{K5}):
\[
\lb{sde}
\m_n=\pi n+\ve_n(p_{c0}-p_{cn}+O(\ve_n)),
\qqq
\ve_n={1\/2\pi n},
\]
\[
\lb{ape}
e_n^\pm=\pi n+\ve_n(p_0\pm |p_n|+O(\ve_n)), \qq \qq
p_n=\int_0^1p(x)e^{-i2\pi nx}dx=p_{cn}-ip_{sn}.
\]

{\bf 2.2. The quasimomentum.} In order to describe the quasimomentum
from \cite{F3}, \cite{MO} we start from the properties of the
Lyapunov function defined by $\D(z)={1\/2}(\vp'(1,z)+\vt(1,z))$. We
shortly describe the properties the Lyapunov function:

1) The function $\D(z)$ is entire and satisfies
(due to the symmetry principle, since $\D$ is real on
$\R$ and $i\R$)
\[
\D(z)=\D(-z)=\ol\D(-\ol z)=\ol\D(\ol z), \qqq\ z\in \C.
\]

2) For each $ n\in \Z$ there exists an unique point $e_n\in [e_n^-,e_n^+]$ such that
\[
\lb{dzn}
\D'(e_n)=0,\qqq (-1)^n\D(e_n)=\max_{\l\in [e_n^-,e_n^+]} |\D(z)|=\cosh h_n\ge 1, \qq \ {\rm for \ some} \ h_n\ge 0.
\]
Recall that $g_n=(e_n^-,e_n^+)$.
Note that if $g_n=\es$, then  $e_n=e_n^\pm$ and  if  $g_n\ne \es$,
then $e_n\in g_n$ and the point $e_n$ is very close to the centrum
of the gap $g_n$ for $n$ large enough, see \er{aen}.

3) $\D(e_n^\pm)=(-1)^n$ and  $\D([e_{n}^+, e_{n+1}^-])=[-1,1]$ for all $n$.

We introduce the quasimomentum $k(\cdot )$ for $H_0$ by
$k(z)=\arccos \D(z),\ z \in \cZ=\C\sm \cup \ol g_n$
and by the asymptotics:
\[
\lb{pk1}
 k(z)=z+O(1/z)\qq  as \ \ |z|\to \iy.
\]
The function $k(z)$ is analytic in $z\in\cZ$.  Moreover, the
quasimomentum $k(\cdot)$ is a conformal mapping from $\cZ$ onto the
quasimomentum domain $\cK$ (see Fig. 5 and \cite{F3}, \cite{MO}) given by
$$
\cK=\C\sm \cup \G_n, \qqq \G_n=(\pi n-ih_n,\pi n+ih_n).
$$
Here $\G_n$ is the vertical cut with the height $h_n$, which is
defined by \er{dzn}. Moreover, we have $(nh_n)_{n\ge 1}\in \ell^2$
iff $p\in L^2(0,1)$, see \cite{MO}, \cite{K1}, \cite{K2}. The
function $\sin k(\cdot)$ is analytic on $\cM$.

We  shortly recall properties of the quasimomentum
$k(\cdot)$ from \cite{MO} or \cite{KK}:

{\bf Properties of the quasimomentum:}

{\it Here and below we rewrite the quasimomentum $k(\cdot)$ in terms
of real functions $u(\cdot), v(\cdot)$ by
$$
k(z)=u(z)+iv(z), \qqq z\in \cZ,
$$
where $u, v$ are harmonic  in $\cZ$. Moreover, $\pm v$ is positive in
$\C_\pm$ and satisfies:
\[
\lb{pfv}
v(z)=\Im z\rt(1+{1\/\pi}\int_{\cup_{n\ne 0}g_n} {v(\t)\/|\t-z|^2}d\t\rt),
\ \ \ z\in \C_\pm.
\]

\no 1)  $v(z)\ge \Im z>0$ and $v(z)=-v(\ol z)$ for all $z\in \C_+$ and
\[
\lb{pk4}
k(z)=-k(-z)=\ol k(\ol z)=-\ol k(-\ol z) ,\qq \forall \ z\in \cZ,
\]
\[
\lb{35} (-1)^{n+1}i\sin k(z)=\sinh v(z)=\pm |\D^2(z)-1|^{1\/2}>0\qq
\ all \qq z\in  g_n^\pm.
\]
\no 2) $v(z)=0$ for all $z\in \s_n=[e_{n-1}^+,e_n^-], n\ge 1$.

\no 3)  If some $g_n\ne \es, n\ge 1$, then the function $v(z+i0)>0$
and $v''(z+i0)<0$ for all $z\in g_n$, and $v(z+i0)$ has  a maximum
at $e_n\in g_n$ (see \er{dzn} and Fig. 4) such that $v'(e_n)=0$, and
\[
\lb{prq}
v(z+i0)=-v(z-i0)>v_n(z)=|(z-e_n^-)(z-e_n^+)|^{1\/2}>0, \qqq
\qqq  \forall \  z\in g_n\ne \es,
\]
\[
\lb{vvn} v(z+i0)=v_n(z)\rt(1+{1\/\pi}\int_{\R\sm
g_n}{v(t+i0)dt\/v_n(t)|t-z|}\rt), \qqq  z\in g_n,
\]
\[
\lb{pk5}
|g_n|\le 2h_n.
\]
\no 4) $u'(z)>0$ on  all $(e_{n-1}^+,e_n^-)$ and  $u(z)=\pi n$ for
all $z\in g_n\ne \es, n\in \Z$.

\no 5) The function $k(z)$ maps a horizontal cut (a "gap" )
$[e_n^-,e_n^+]$ onto the vertical cut $\G_n$  and the momentum band
$\s_n=[e_{n-1}^+,e_n^-]$ onto the segment $[\pi (n-1), \pi n]$ for
all $n\in\Z$, i.e.,
\[
\lb{pk6} k([e_n^-,e_n^+])=\G_n, \qqq k(\s_n)=[\pi (n-1), \pi n], \qqq n\in \Z.
\]

\no 6) The following  identities hold true:
\[
\lb{pk7}
k(z)=z+{1\/\pi}\int_{\cup g_n} {v(t+i0)dt\/t-z}, \qqq z\in \ol\C_+\sm \cup \ol g_n.
\]

}

{\bf 2.3. The momentum Riemann surface $\cM$}. Recall that we will
work with the momentum $z=\sqrt \l$, where $\l\in \L$ is an energy.
The function $\l\to z=\sqrt \l$ maps the  cut Riemann surface $\L\sm
\cup \g_n^c$ onto  the cut  momentum domain $\cZ$ given by
 \[
\lb{Z00}
\cZ=\C\sm \cup_{n\ne 0} \ol g_n, \qq g_n=(e_n^-,e_n^+)=-g_{-n},
\qq e_n^\pm=-e_{-n}^\mp=\sqrt{E_n^\pm}>0,\qq
n\ge 1.
\]
Slitting the n-th nontrivial momentum gap  $g_n$, we obtain a cut
$g_n^c$ with an upper   $g_n^+$ and lower rim  $g_n^-$. Below we
will identify this cut $g_n^c$ and the union of
  the upper rim  (gap) $\ol g_{n}^+$ and the lower rim (gap)
  $\ol g_{n}^{\ -}$, i.e.,
\[
\lb{gn}
g_n^c=\ol g_{n}^+\cup \ol g_{n}^-,\ \ where \ \ g_{n}^\pm =g_n\pm i0.
\]
In order to construct the  Riemann surface $\cM$ we take the cut
domain $\cZ$ and identify (i.e. we glue) the upper rim $g_{n}^+$ of
the cut $g_n^c$ with the upper rim $g_{-n}^+$  of the cut $g_{-n}^c$
and correspondingly the lower rim $g_{n}^-$ of the cut $g_n^c$ with
the  lower  rim $g_{-n}^-$  of the cut $g_{-n}^c$ for all nontrivial
gaps. The mapping $z=\sqrt \l: \L\to \cM$ is one-to-one and onto.
The bounded physical  gap $\g_n^{(1)}\ss \L_1$ is mapped onto
$g_n^+\ss \cM_1$ and the bounded nonphysical gap $\g_n^{(2)}\ss
\L_2$ is mapped onto $g_n^-\ss \cM_2$. Moreover,

1) $\C _+$ plus all physical gaps $g_{n}^+$ is a so-called physical
"sheet" $\cM_1$,

2) $\C _-$ plus all nonphysical gaps $g_{n}^-$ is a so-called
non physical "sheet" $\cM_2$.

3) The  momentum bands $\s_n=[e_{n-1}^+,e_n^-], n\in \Z$   joint
the first and second sheets.

For the construction of the  Riemann surface $\cM$ we need
to write few simple remarks:

1) $\cM$ is the  Riemann surface of the function $\sin k(z)=\sqrt{1-\D^2(z)}$
and $\L$ is the  Riemann surface of the function
$\sin k(\sqrt \l)=\sqrt{1-\D^2(\sqrt \l)}$.

2) It is important that \er{pk4} gives
$$
\sin k(z+i0)=\sin k(-z+i0)=\sin (\pi n+iv)
=(-1)^ni\sinh v,\ \ v=\Im k(z)>0, \ \ \forall \ z\in g_n\ne \es.
$$
Due to this identity the upper rim $g_{n}^+$ of the cut $g_n^c$
is glued with the upper rim $g_{-n}^+$  of the cut $g_{-n}^c$.
Correspondingly the lower rim $g_{n}^-$ of the cut $g_n^c$ is glued
with the  lower  rim $g_{-n}^-$  of the cut $g_{-n}^c$.
Due to these facts the function $\sin k(z), z\in \cZ$  is analytic
in $\cM$ and $\cM$ is the Riemann surface of $\sin k(z)$.

3)  Let $f$ be entire. The function $f(z), z\in \C$ is even,
i.e., $f(z)=f(-z), z\in \C$, iff $f$ is analytic on
the  Riemann surface $\cM\ne \C$.

\begin{figure}
\unitlength 1mm 
\linethickness{0.4pt}
\ifx\plotpoint\undefined\newsavebox{\plotpoint}\fi 
\begin{picture}(119.75,66.5)(0,0)
\put(15,17.25){\line(1,0){104.75}}
\qbezier(40.5,17.25)(40.5,30.375)(52.25,35.75)
\qbezier(52.25,35.75)(58.25,38.375)(64.25,35.75)
\qbezier(64.25,35.75)(76.17,30.375)(76,17.25)
\qbezier(40.5,17.25)(40.5,45.125)(58.25,46.5)
\qbezier(76,17.25)(76.17,45.125)(58.25,46.5)
\put(38.75,12.00){\makebox(0,0)[cc]{$z_n^-$}}
\put(76,12){\makebox(0,0)[cc]{$z_n^+$}}
\put(65.25,38.5){\makebox(0,0)[cc]{$v_n$}}
\put(70,44.25){\makebox(0,0)[cc]{$v$}}
\put(18,17.25){\line(0,-1){1.00}}
\put(95,17.25){\line(0,-1){1.00}}
\put(57,17.25){\line(0,-1){1.00}}
\put(18,12.00){\makebox(0,0)[cc]{$z_{n-1}^+$}}
\put(95,12.00){\makebox(0,0)[cc]{$z_{n+1}^-$}}
\put(57,12.00){\makebox(0,0)[cc]{$e_{n}$}}
\end{picture}
\caption{The graph of $v(z+i0), \ z\in [z_{n-1}^+, z_{n+1}^-]$
and $h_n=v(e_n+i0)>0$}
\lb{grafv}
\end{figure}
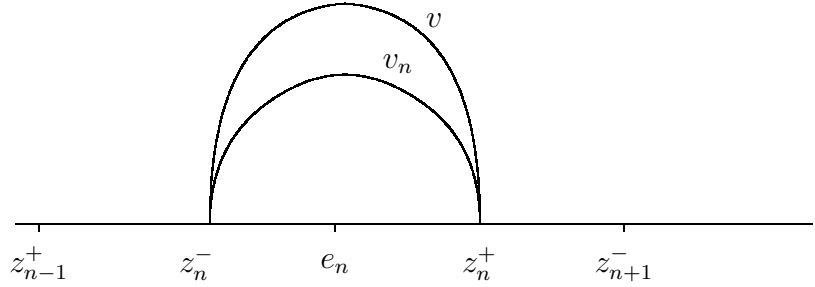

{\bf 2.4. The Floquet solutions.}
The Floquet solutions $\p_{\pm}(x,z), z \in \cZ$ of $H_0$ are given by
\[
\lb{3} \p_\pm(x,z)=\vt(x,z)+m_\pm(z)\vp(x,z),\qq m_\pm={\b\pm i\sin
k\/ \vp(1,\cdot)},\qqq \b={\vp'(1,\cdot)-\vt(1,\cdot)\/2},
\]
where $\vp(1,z)\p_+(\cdot,z)\in L^2(\R_+)$ for all
$z\in\C_+\cup\cup_{n\ne 0} g_n^+$. If $p=0$, then $k=z$ and
$\p_\pm(x,z)=e^{\pm izx}$. Substituting estimates \er{efs1} into
\er{fs}  we obtain the standard asymptotics
\begin{multline}
\lb{asb}
\b(z)=\int_0^1{\sin z(2x-1)\/z}p(x)dx+{O(e^{|\Im
z|})\/z^2},\\ \dot\b(z)=\int_0^1{\cos z(2x-1)\/z} p(x)(2x-1)
dx+{O(e^{|\Im z|})\/z^2}\qq as \qq |z|\to \iy, \ here \ \ \dot \b=\pa_z \b.
\end{multline}
\begin{figure}
\tiny
\unitlength=1mm
\special{em:linewidth 0.4pt}
\linethickness{0.4pt}
\begin{picture}(120.67,34.33)
\put(20.33,20.00){\line(1,0){102.33}}
\put(71.00,7.00){\line(0,1){27.00}}
\put(70.00,18.67){\makebox(0,0)[cc]{$0$}}
\put(124.00,18.00){\makebox(0,0)[cc]{$\Re k$}}
\put(67.00,33.67){\makebox(0,0)[cc]{$\Im k$}}
\put(87.00,15.00){\linethickness{2.0pt}\line(0,1){10.}}
\put(103.00,17.00){\linethickness{2.0pt}\line(0,1){6.}}
\put(119.00,18.00){\linethickness{2.0pt}\line(0,1){4.}}
\put(56.00,15.00){\linethickness{2.0pt}\line(0,1){10.}}
\put(39.00,17.00){\linethickness{2.0pt}\line(0,1){6.}}
\put(23.00,18.00){\linethickness{2.0pt}\line(0,1){4.}}
\put(85.50,18.50){\makebox(0,0)[cc]{$\pi$}}
\put(54.00,18.50){\makebox(0,0)[cc]{$-\pi$}}
\put(101.00,18.50){\makebox(0,0)[cc]{$2\pi$}}
\put(36.00,18.50){\makebox(0,0)[cc]{$-2\pi$}}
\put(117.00,18.50){\makebox(0,0)[cc]{$3\pi$}}
\put(20.00,18.50){\makebox(0,0)[cc]{$-3\pi$}}
\put(87.00,26.00){\makebox(0,0)[cc]{$\pi+ih_1$}}
\put(56.00,26.00){\makebox(0,0)[cc]{$-\pi+ih_1$}}
\put(103.00,24.00){\makebox(0,0)[cc]{$2\pi+ih_2$}}
\put(39.00,24.00){\makebox(0,0)[cc]{$-2\pi+ih_2$}}
\put(119.00,23.00){\makebox(0,0)[cc]{$3\pi+ih_3$}}
\put(23.00,23.00){\makebox(0,0)[cc]{$-3\pi+ih_3$}}
\end{picture}
\caption{The domain $\cK=\C\sm \cup \G_n$ with
the cuts $\G_n=(\pi n-ih_n,\pi n+ih_n)$}
\lb{k}
\end{figure}
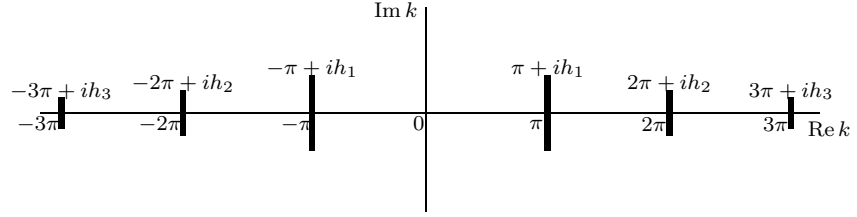
The function $\sin k$ and each function $\vp(1,\cdot)\p_{\pm}(x,\cdot),
x\in \R$ are analytic on the Riemann surface $\cM$. Recall that
the Floquet solutions $\p_\pm(x,z), (x,z)\in \R\ts \cM$ satisfy
(see \cite{T})
\[
\lb{f1}
\p_\pm(0,z)=1, \qq \p_\pm(0,z)'=m_\pm(z),
\qq \p_\pm(1,z)=e^{\pm ik(z)}, \qq  \p_\pm(1,z)'=e^{\pm ik(z)}m_\pm(z),
\]
\[
\lb{f2}
\p_\pm(x,z)=e^{\pm ik(z)x}(1+O(1/z)) \qq as \qq |z|\to \iy, \qq z\in \cZ_\ve,
\]
uniformly in $x\in \R$,  where the set $\cZ_\ve$ is given by
$$
\cZ_\ve =\{z\in \cZ: \dist \{z,g_n\}>\ve, g_n\ne \es,  n\in \Z\},\qq \ve>0.
$$
Below we need the simple identities
\[
\lb{LD0}
\b^2+1-\D^2=1-\vp'(1,\cdot)\vt(1,\cdot)= -\vp(1,\cdot)\vt'(1,\cdot).
\]
This yields
\[
\lb{wtf}
m_+(z)m_-(z)= -{\vt'(1,z)\/\vp(1,z)},\qqq z\ne \m_n.
\]
Let $\dD_r(z_0)=\{|z-z_0|<r\}$ be a disk for some $r>0$.
We need the following result (see \cite{Zh3}), where
$\mA(z_0), z_0\in \cM$, denotes the set of functions analytic in some
disc $\dD_r(z_0), r>0$.

\begin{lemma}
\lb{Tm}
i) The following asymptotics hold true:
\[
\lb{Tm-1}
m_\pm (z)=\pm iz+O(1) \qq as \qq |z|\to \iy, \qq z\in \cZ_\ve,\ve >0.
\]

\no ii) If $g_n=\es$ for some $n\in \Z$, then
the functions $\sin k, m_\pm$ are analytic in some  disk
$\dD(\m_n,\ve)\ss\cZ, \ve>0$. The functions $\sin k(z)$
and $\vp(1,z)$ have the simple zero at $\m_n$ and satisfies
\[
\lb{Tm-2}
m_\pm (\m_n)={\dot\b(\m_n)\pm i(-1)^n\dot k(\m_n)\/\dot\vp(1,\m_n)},
\qq \Im m_\pm (\m_n)\ne 0.
\]

\no iii) If the function $m_+$ has a pole at $\m_n+i0$
 for some $n\ge 1$, then  $k(\m_n+i0)=\pi n+ih_{sn}$ and
\begin{multline}
\lb{Tm-31}
\qq h_{sn}>0,\qq
\b(\m_n)=i\sin k(\m_n+i0)= -(-1)^n\sinh h_{sn},\qq m_+\in \mA(\m_n-i0),\\
 \qq m_+(\m_n+z)={\r_n+O(z)\/z} \qq as \ z\to 0, \ z\in \C_+,\
\qq \r_n={-2\sinh |h_{sn}|\/(-1)^n\dot\vp(1,\m_n)}<0.
\end{multline}

\no iv)  If the function $m_+$ has a pole at $\m_n-i0$
for some $n\ge 1$, then $k(\m_n-i0)=\pi n+ih_{sn}$ and
\begin{multline}
\lb{Tm-32}
h_{sn}<0,\qq
\b(\m_n)=-i\sin k(\m_n-i0)=(-1)^n\sinh h_{sn},\ m_+\in \mA(\m_n+i0),\\
 \qq m_+(\m_n+z)={\r_n+O(z)\/z}\qq as \ z\to 0, \ z\in \C_-.
\end{multline}

\no v)  Let $e_n^-< e_n^+$
for some $n\ne 0$. Then   $\m_n=e_n^-$ or $\m_n=e_n^+$  iff
\[
m_+(\m_n+z)={\r_n^\pm+O(z)\/\sqrt z}\qq as \ z\to 0, \
 z\in \C_+,\qq\ {\rm for \ some} \ \const \
\r_n^\pm\ne 0.
\]
\end{lemma}

\no {\bf Proof.} The results of this lemma is well-known.
We will give a sketch.

i) The asymptotics \lb{Tm-1} follows from \er{asb}, \er{efs1},
\er{pk1},  see \cite{T}.

ii) If $g_n=(e_n^-,e_n^+)=\es$, then due to \er{pk6}  we have $k(e_n^\pm)=\pi n$ and the function $k(\cdot)$ is analytic at $e_n^-=e_n^+$. Moreover, \er{pk7} gives
$$k'(z)=1+{1\/\pi}\int_{\cup g_n} {v(t+i0)dt\/(t-z)^2}\ge 1,\qqq {\rm at}
\qqq z=e_n^\pm.
$$  Thus this yields that the function $\sin k(z)$ is analytic at $z=e_n^\pm$ and $z=e_n^\pm$ is a simple
zero of $\sin k(z)$. The point $\m_n=e_n^\pm$ is a simple zero of  $\vp(1,z)$, since the point $\m_n^2$ is the Dirichlet eigenvalue, see \er{Dei}. This implies the proof of ii).

iii) Let  $g_n=(e_n^-,e_n^+)\ne \es$ and let the function $m_+$ have a pole at $\m_n+i0$.
The point $\m_n\in g_n$, then $k(\m_n+i0)=\pi n+ih_{sn}\in \G_n$ for some $h_{sn}>0$, since
$k$ is the conformal mapping and $k(g_n)=\G_n$ and $\m_n\ne e_n^\pm$.

Due to \er{wtf}, the function $m_-$ is analytic at $\m_n+i0$ and we get
$m_-(\m_n+i0)\ne 0$. Then $\b(\m_n)-i\sin k(\m_n+i0)=0$, which yields
$$
\b(\m_n)=i\sin k(\m_n+i0)=i\sin (\pi n+ih_{sn})={(-1)^n\/2}(e^{-h_{sn}}-e^{h_{sn}})=-(-1)^n\sinh h_{sn}.
$$
Moreover, using similar arguments ( analytic properties of $k(\cdot)$), we obtain
$m_+(\m-i0)=m_-(\m+i0)$ for all $\m\in (e_n^-,e_n^+)$, this gives $m_+\in \mA(\m_n-i0)$.

Let $ z\to 0, \ z\in \C_+ $. Then
$$
m_+(\m_n+z)={\b(\m_n+z)+i\sin k(\m_n+z)\/\vp(1, \m_n+z)}={-2(-1)^n\sinh h_{sn}+O(z)\/z(\dot\vp(1, \m_n)+O(z))}=
{\r_n+O(z)\/z}.
$$
Thus iii) has been proved.

The proof of iv) and v) is similar to the case of iii).
\BBox

{\bf 2.5. Properties of fundamental solutions.} Let $\n_n^2, n\ge 1$
be the Neumann  spectrum of the equation $-y''+py=\n^2y$ on the
interval $[0,1]$ with the boundary condition $y'(0)=y'(1)=0$. It is
well known that each $\n_n^2\in [E^-_n,E^+_n ], n\ge 1$.

\begin{lemma}
\lb{T31}
Let $p\in L^1(0,1)$. Then

\no i) The following asymptotics hold true uniformly in
$z\in [e_n^-,e_n^+]$ as $n\to \iy$:
\[
\lb{T31-1}
\vp(1,z)=(-1)^n{(z-\m_n)\/\pi n}(1+O(1/n)),\qq
\]
\[
\lb{T31-2}
-{\vt'(1,z)\/z^2}=(-1)^n{(z-\n_n)\/\pi n}(1+O(1/n)).
\]
 ii) Let $z\in g_n\ss (8e^{\|p\|_1},\iy)$.  Then the following
 estimates hold true (here $\dot y=\pa_z y$):
\[
\lb{T31-3}  |g_n|^2\le {8e^{\|p\|_1}\/|z_n|}<1,
\]
\[
\lb{T31-4}
|\dot \vp(1,z)|\le {3e^{\|p\|_1}\/2|z|},
\qqq |\vp(1,z)|_{z\in g_n} \le |g_n|{3e^{\|p\|_1}\/2|z|},\qq
\]
\[
\lb{T31-5}
|\dot \vt'(1,z)|\le |z|{3e^{\|p\|_1}\/2},\qqq
|\vt'(1,z)|\le |g_n||z|{3e^{\|p\|_1}\/2},
\]
\[
\lb{T31-6}
|\dot\b(z)|\le  {3e^{\|p\|_1}\/2|z|}, \qqq
|\b(e_n^\pm)|\le |g_n|{3e^{\|p\|_1}\/2},\qqq
|\b(z)|\le  |g_n|{9e^{\|p\|_1}\/4|z|}.
\]
\no iii) In each disk $\dD_{\pi\/4}(\pi n)\ss \dD=\{|z|>32e^{2\|p\|_1}\}$
there exists exactly one momentum gap $g_n$. Moreover, if
$g_n,g_{n+1}\ss\dD$, then $e_{n+1}^--e_{n}^+\ge \pi$.
\end{lemma}
\no {\bf Proof.} i) We have the Taylor formula
$\vp(1,z)=\dot\vp(1,\m_n)\t+\ddot\vp(1,\m_n+\a \t){\t^2\/2}$
for any $z\in [e_n^-, e_n^+]$ and some $\a\in [0,1]$, where $\t=z-\m_n$.
Asymptotics \er{efs1} give $\dot\vp(1,\m_n)=2(-1)^n{(1+O(1/n))\/2\pi n}$
and $\ddot\vp(1,\m_n+\a \t)\t=O(n^{-2})$, which yields \er{T31-1}.
Similar arguments imply \er{T31-2}.

ii) Using $|\D(z)-\cos z|\le {e^{\|p\|_1}\/|z|}$, for all $|z|\ge 2$, we obtain
$$
{h_n^2\/2}\le \cosh h_n-1=|\D(z_n)|-1\le {e^{\|p\|_1}\/|z_n|}.
$$
Then the estimate $|g_n|\le 2h_n$ (see \er{pk5}) gives  \er{T31-3}.

Due to \er{efs1}, the function $f=z\vp(1,z)$ has the estimate
$|f(z)|\le C_1=e^{\|p\|_1}, z\in \R$. Then
the Bernstein inequality gives $|\dot f(z)|=|\vp(1,z)+
z\dot \vp(1,z)|\le C_1, z\in \R$, which yields
$|\dot \vp(1,z)|\le {C_1\/|z|}(1+{1\/|z|}), z\in \R$.
Moreover, we obtain $|\vp(1,z)|\le |g_n|\max_{z\in g_n}
 |\dot \vp(1,z)|\le |g_n|{3C_1\/2|z|}$.

The proof of \er{T31-5} and the estimate $|\dot\b(z)|\le
 {3e^{\|p\|_1}\/2|z|}$ is similar. Identity \er{LD0} gives

$\b^2(e_n^\pm)=-\vp(1,e_n^\pm)\vt'(1,e_n^\pm)$.
Then \er{T31-4}, \er{T31-5} imply $|\b(e_n^\pm)|\le |g_n|{3e^{\|p\|_1}\/2}$.
Using these estimates and $\b(z)=\b(e_n^-)+\dot\b(z_*)(z-e_n^-)$ for some
$z_*\in g_n$, we obtain \er{T31-6}.

iii) Using \er{efs1} we obtain
$$
|(\D^2(z)-1)-(\cos^2z-1)|\le 2X|\D(z)-\cos z|\le 2X^2/|z|,\qq
X=e^{|\Im |+\|p\|_1}.
$$
After this the standard arguments (due to Rouche's theorem) give
the proof of iii).
\BBox

\begin{lemma}
\lb{Ted} i) Let $q\in L^1(0,t)$ for some $t>0$ and let
$z-e_n^\pm=O(|g_n|/n)$ as $n\to \iy$. Then
\[
\lb{dex}
\int_\R\vp(1,z,\t)q(\t)d\t={(-1)^n|g_n|\/2\pi n}
(\mp \wh q_0+c_n\wh q_{cn}-s_n\wh q_{sn}+ O(1/n)), \qq
\]
where $c_n, s_n$ are defined in \er{cs} and $\wh q_0=\int_\R q(\t)d\t$.

ii) Let  $v_n$ be defined by \er{prq}.  The following asymptotics holds true:
\[
\lb{aen}
e_n={e_n^-+e_n^+\/2}+O(|g_n|^2/n^2),
\]
\[
\lb{pav}
v(z)=\Im k(z)=\pm v_n(z)(1+O(n^{-2})),\qqq   z\in g_n\pm i0\qq as \ n\to \iy.
\]
\end{lemma}
\no {\bf Proof.} i) We consider the case $z-e_n^-=O(|g_n|/n)$,
the proof for $z-e_n^+=O(|g_n|/n)$
  is similar.

  We need the following facts from Theorem 2 in \cite{K5}:
Let $\m_n^2(\t), \t\in \R$ be the Dirichlet eigenvalue for
the problem $-y''+q(x+\t)y=z^2 y, y(0)=y(1)=0$.
Then there exists a real function  $\f_n(\t), \t\in \R$  such
that $\f_n',\f_n''\in L_{loc}^2(\R)$ and the following identity and asymptotics
\[
\lb{mz1}
{E_n^-+E_n^+\/2}-\m_n^2(\t)={|\g_n|\/2}\cos \f_n(\t),\qq \qq
\forall \ \t\in \R,\qq \cos \f_n(0)=c_n, \qq
\sin \f_n(0)=s_n,
\]
\[
\lb{mz2} \f_n(\t)=\f_n(0)+2\pi n\t+O(\ve_n)\qqq as \qqq n\to
\iy,\qqq \ve_n={1\/2\pi n},
\]
hold true, uniformly with respect to $\t\in [0,1]$.

Using \er{mz1} we rewrite $e_n^--\m_n(\t)$ in the form
\begin{multline}
\lb{a2}
e_n^--\m_n(\t)={{E_n^-+E_n^+-|\g_n|\/2}-\m_n^2(\t)\/e_n^-+\m_n(\t)}
=
{|\g_n|\/2} {(-1+\cos \f_n(\t))\/e_n^-+\m_n(\t)}
\\
={|g_n|\/2}(-1+\cos \f_n(\t)+O(|g_n|\ve_n)),
\end{multline}
where the following asymptotics  have been used:
\[
\lb{a3} {|\g_n|\/e_n^-+\m_n(\t)}=|g_n|{e_n^-+e_n^+\/e_n^-+\m_n(\t)}
=|g_n|\rt(1+  {e_n^+-\m_n(\t)\/e_n^-+\m_n(\t)}\rt) =|g_n|(1+
O(|g_n|\ve_n))
\]
as $n\to \iy$.
Asymptotics \er{T31-1} and \er{a2} yield
\[
\lb{a1} \vp(1,z,\t)={(-1)^n\/\pi n}(1+O(\ve_n))(z-\m_n(\t))
={(-1)^n\/\pi n}(1+O(\ve))(e_n^--\m_n(\t)+\ve_n O(|g_n|))
\]
$$
={(-1)^n|g_n|\/2\pi n}(-1+\cos \f_n(\t)+O(\ve_n)).
$$
Combine \er{a1}, \er{mz2} we obtain
$$
\int_\R\vp(1,z,\t)q(\t)d\t ={(-1)^n|g_n|\/2\pi n} \int_\R
\rt[-1+\cos \f_n(\t)+ O(\ve_n)\rt]q(\t)d\t
$$
$$
={(-1)^n|g_n|\/2\pi n}(-\wh q_0+c_n\wh q_{cn}-s_n\wh q_{sn}+
O(\ve_n)),
$$
where $c_n=\cos \f_n(0), s_n=\sin \f_n(0)$, and this yields \er{dex}.

ii) We need the estimate from Lemma 2.1 in \cite{K4} for all $n\ge
1$:
\begin{multline}
\lb{een}
  \rt|e_n-{e_n^-+e_{n}^+\/2}\rt|\le {|g_n|^2\/4}M_n,\qqq
  M_n=\max _{z\in \ol g_n}G_n(z),\\
\qqq  G_n(z)={1\/\pi}\int_{\R\sm g_n}{v(\t+i0)d\t\/v_n(\t)|\t-z|},\
\ z\in g_n.
\end{multline}
We need the estimate from Lemma 2.4 in \cite{K4}:
\[
\lb{eGn}
 G_n(z)\le {1\/|z|\sqrt{|e_n^+e_n^-|}}\rt({\wh p_0}
 +{\int_0^1p^2(x)dx\/4r^2}\rt),\qqq  \ z\in \ol g_n,
 \qqq r=\min_{n} |e_n^+-e_{n+1}^-|>0.
\]
This gives $M_n=O(1/n^2)$ as $n\to \iy$, since
$e_n^\pm =\pi n+o(1)$ as $n\to \iy$.
Thus substituting \er{eGn} into \er{een} and into \er{vvn},
we obtain  \er{aen} and  \er{pav}.
\BBox

\section {The perturbed operator $H$ }
\setcounter{equation}{0}

We recall the well-known results about the scattering for $H, H_0$,
see e.g. \cite{F3}, \cite{F1}. The equation $-f''+(p+q)f=z^2 f$ has
unique Jost solutions $f_\pm (x,z)$ such that
\[
\lb{bcf}
f_+(x,z)=\p_+(x,z), \ x\ge t, \ \ {\rm and} \ \
f_-(x,z)=\p_-(x,z), \ x\le 0, \qq z\in \s_M=\R\sm \cup [e_n^-,e_n^+].
\]
The Jost solutions satisfy
 $$
 f_+(x,z)=\ol f_+(x,-z), \qqq \qqq  \forall z\in\s_M.
 $$
  This yields the basic identity
\[
f_+(x,z)=b(z)f_-(x,z)+a(z)f_-(x,-z),\qqq  \forall z\in\s_M,
\]
 where
\begin{multline}
\lb{abws}
b={s\/w_0}, \qqq a={w\/w_0},\qqq s=\{f_+(x,z),f_-(x,-z)\},\\
\ w=\{f_-,f_+\},  \qqq w_0=\{\p_-,\p_+\}={2i\sin k\/\vp(1,\cdot)},\qqq\qqq
\end{multline}
and $\{f, g\}=fg'-f'g$ is the Wronskian.
The scattering matrix $\cS_M$ for $H, H_0$ is given by
\[
\cS_M (z)\ev \ma a(z)^{-1}& r_-(z)\\ r_+(z)&a(z)^{-1}\am ,
\ \ \ \  r_{\pm}={s(\mp z)\/ w(z)}=\mp {b(\mp z)\/ a(z)}, \ \ \ \  z\in \s_M,
\lb{1.7}
\]
where $1/a$ is the transmission coefficient and $r_{\pm}$
is the reflection coefficient. We have the following identities
from \cite{F1}, \cite{F3}:
\[
\lb{iab}
|a(z)|^2=1+|b(z)|^2,\qqq z\in \s_M=\R\sm \cup [e_n^-,e_n^+].
\]

We will represent the Jost solutions $f_\pm (x,z)$ in the form
$f_+=\wt\vt +m_+\wt\vp$ (see \er{T22-1}) and recall that $m_\pm$ is
the Weyl-Titchmarsh function given by \er{3}. Here $\wt\vt, \wt\vp$
are the solutions of the equations $-y''+(p+q)y=z^2y, z\in \C$ and
satisfying
\[
\lb{wtc}
\wt\vp(x,z)=\vp(x,z),\qqq \wt\vt(x,z)=\vt(x,z) \ \ for \ \ all \ x\ge t.
\]
Due to \er{ieq},  the solutions $\wt\vt, \wt\vp$ and $f_+$  of the equation
$-y''+(p+q)y=z^2y$ satisfy the equation
\[
\lb{ep}
y(x,z)=y_0(x,z)-\int_x^t\vp(x-\t,z,\t)q(\t)y(\t,z)d\t, \qqq x\le t,
\]
where $y$ is one from  $\wt\vt, \wt\vp$ and $f_+$; $y_0$ is the
corresponding function from  $\vt, \vp$ and $\p_+$. For each $x\in
\R$ the functions  $\wt\vt (x,z), \wt\vp (x,z)$  are entire in
$z\in\C$ and satisfy
\begin{multline}
\lb{efs}
\max \rt\{|z|_1|\wt\vp(x,z)|, \ |\wt\vp'(x,z)| , |\wt\vt(x,z)|,
{|\wt\vt'(x,z)|\/|z|_1}    \rt\} \le X_1,\\
|\wt\vt(x,z)-\vt(x,z)|\le {\|q\|_t\/|z|_1}X_1,\qq
|\wt\vp(x,z)-\vp(x,z)|\le {\|q\|_t\/|z|_1^2}X_1,\\
\|p\|_t=\int_0^t|p(s)|ds,\qq
|z|_1=\max\{1, |z|\},\qqq X_1=e^{|\Im z||2t-x|+\|q\|_t+\|p\|_t+\int_x^t|p(\t)|d\t},
\end{multline}
for all $(p,x,z)\in L_{loc}^1(\R)\ts \R\ts \C$. The proof of
\er{efs} repeats the standard arguments (see \cite{PT})
proving \er{efs1}.

\begin{lemma}
\lb{T22} Each function  $f_+(x,z), x\in \R$ has an analytic
continuation from $z\in \s_M$ into $z\in\cZ$. Moreover, for all
$z\in \cZ$ the following identities and asymptotics  hold true:
\[
\lb{T22-1}
f_+(\cdot,z)=\wt\vt (\cdot,z)+m_+(z)\wt\vp (\cdot,z),
\]
\begin{multline}
\lb{T22-2}
f_+(0,z)\!\!=1\!+\!\int_0^t\!\vp(x,z)q(x)f_+(x,z)dx,\qq
f_+'(0,z)\!\!=\! m_+(z)-\!\int_0^t\! \vt(x,z)q(x)f_+(x,z)dx,
\end{multline}
\[
\lb{T22-20}
|f_+(x,z)-\p_+(x,z)|\le e^{-vx+B(t,x)}{\e(z)\/|z|_1}\int_x^t|q(r)|dr,\qq
\forall x\in [0,t],
\]
where $v=\Im k(z)$ and
$$
\e(z)=\sup_{x\in [0,t]} |e^{-ik(z)x}\p_+(x,z)|,
\qq B(t,x)=2(t-x)\gJ+\int_x^t (|p(r)|+|q(r)|)dr, \qq \gJ={|v|-v\/2},
$$
\[
\lb{T22-3} f_+(x,z)=e^{\pm ik(z)x}(1+e^{\pm (t-x)2\gJ}O(1/z)), \qq
    \qq x\in [0,t],
\]
\begin{multline}
\lb{T22-4}
\e(z)\to 1,\qqq
f_+(0,z)=1+{O(e^{2t\gJ})\/z},   \qqq f_+'(0,z)=iz+O(1)+o(e^{2t\gJ})
\end{multline}
as $|z|\to \iy,\ \ z\in \cZ_\ve,\ \ve>0$, where
$\wh q(z)=\int_0^tq(x)e^{2izx}dx, \ z\in \C$.

\end{lemma}
\no{\bf Proof.} Using \er{wtc}, \er{bcf} we obtain \er{T22-1}.
Then each function  $f_+(x,z), x\in \R$  is analytic in $z\in\cZ$, since
$m_+$ is analytic in $z\in\cZ$.

Using the identity \er{ffs}, we obtain
$\vp(-s,\cdot,s)=-\vp(s,\cdot)$ and $\vp'(-s,\cdot,s)=\vt(s,\cdot)$.
Substituting the last identities into \er{ep} we get \er{T22-2}.

We will show \er{T22-20} for the case $z\in \ol\C_-$,
the proof for $z\in \C_+$ is similar. Let
$$
f(x)=e^{-ik(z)x}f_+(x,z),\qq f_0(x)=e^{-ik(z)x}\p_+(x,z),
\qq K(x,s)=\vp(x-s,z,s)e^{ik(z)(s-x)}.
$$
The standard iteration of \er{ep} yields
\[
\lb{eU}
f=f_0+\sum _{n\ge 1}f_n, \qqq f_n(x,z)=-\int_x^tK(x,s)q(s)f_{n-1}(s,z)ds.
\]
Using \er{efs1} and estimate $|\Im z|\le |v(z)|, z\in \C$
(see \er{pfv}) we obtain
\begin{multline}
\lb{eK}
|e^{ik(z)(s-x)}|\le e^{|v|(s-x)},\qqq |\vp(x-s,z,s)|
\le {e^{|\Im z|(s-x)+\int_x^s |p(r)|dr}\/|z|_1}\le
{e^{|v|(s-x)+\int_x^s |p(r)|dr}\/|z|_1},\\
|K(x,s)|\le {e^{2|v|(s-x)+\int_x^s|p(r)|dr}\/|z|_1},\qqq \qqq
\qqq \forall s>x, z\in \ol\C_-.\qqq \qqq
\end{multline}
  Substituting \er{eK} into \er{eU}  we obtain
$$
|f_n(x)|\le \int_x^t|K(x,s_1)q(s_1)f_{n-1}(s_1)|ds_1 \le
$$
$$
\le \e \int_x^t|K(x,s_1)q(s_1)|ds_1
\int_{s_1}^t|K(s_1,s_2)q(s_2)|ds_2....
\int_{s_{n-1}}^t|K(s_{n-1},s_{n})q(s_{n})|ds_{n}
$$
$$
\le {\e e^{\int_s^t |p(r)|dr}\/|z|_1^n}
\int_x^te^{2|v(z)|(s_1-x)+\int_x^{s_1}|p|dr}|q(s_1)|ds_1
\int_{s_1}^te^{2|v(z)|(s_2-s_1)+\int_{s_1}^{s_2}|p(r)|dr}|q(s_2)|ds_2....
$$$$
\int_{s_{n-1}}^te^{2|v(z)|(s_n-s_{n-1})
+\int_{s_{n-1}}^{s_n}|p(r)|dr}|q(s_{n})|ds_{n}
$$
$$
\le {\e e^{2|v(z)|(t-x)+\int_s^t |p(r)|dr}\/|z|_1^n}
\int_x^t|q(s_1)|ds_1  \int_{s_1}^t|q(s_2)|ds_2....
\int_{s_{n-1}}^t|q(s_{n})|ds_{n}
$$
$$
=\e (z)e^{2|v(z)|(t-x)+\int_s^t |p(r)|dr}{(\int_x^t|q(r)|dr)^n\/n!|z|_1^n}.
$$
 Thus summing we deduce that
$$
|f(x,z)-f_0(x,z)|\le \e e^{2(t-x)\gJ+\int_x^t |p(r)|dr}
\sum_{n\ge 1}{(\int_x^t|q(r)|dr)^n\/|z|_1^nn!}
\le {\e \/|z|_1}e^{B(t,x)}\int_x^t|q(r)|dr,
$$
which yields \er{T22-20}. Substituting \er{f2} into \er{T22-20} we
obtain \er{T22-3}. The proof of \er{T22-4} is similar. \BBox

Firsova  \cite{F4}, \cite{F1} obtained the following results:
the function  $a(z)$ has an analytic continuations
 from $z\in \s_M$ into $z\in\C_+$ and the following identity holds true:
\[
a(z)=D(z)=\det (I+q(H-z^2)^{-1}), \qqq z\in\C_+.
\]
We prove the main result of this section.

\begin{lemma}
\lb{T23}
i) The functions  $\x(z), s(z), w(z), a(z)$ have analytic continuations
 from $z\in \s_M$ into $z\in\cZ$ and satisfy
 \[
\lb{symab}
a(-z)=\ol a(\ol z), \qq w(-z)=\ol w(\ol z),
\qq s(-z)=\ol s(\ol z),\qq w_0(-z)=\ol w_0(\ol z),\qq  \forall \ z\in \cZ.
\]
Moreover, for each $z\in \cZ$ the following identities
\[
\lb{T23-1}
 w(z)=f_+'(0,z)-m_-(z)f_+(0,z)={2i\sin k(z)\/\vp(1,z)}-
 \int_0^tq(x)\p_-(x,z)f_+(x,z)dx,
\]
\[
\lb{iaw}
\x(z)=2ia(z)\sin k(z) =\vp(1,z) w(z),
\]
\[
\lb{T23-2}
s(z)=m_+(z)f_+(0,z)-f_+'(0,z)=\int_0^tq(x)\p_+(x,z)f_+(x,z)dx,
\]
and the following asymptotics
\[
\lb{T23-3}
\x(z)=2i\sin z(1+O(e^{2t\gJ}/z)) ,\qqq s(z)=O(e^{2t\gJ}),
\]
hold true as $|z|\to \iy, z\in \cZ_\ve, \ve>0$.

\no ii) The function $s(\cdot)$ has exponential type $\r_\pm$
in the half plane $\C_\pm$, where $\r_+=0, \r_-=2t$.

\end{lemma}
\no{\bf Proof.} We have $w=\{f_-,f_+\}=\p_-{f_+}'-m_-f_+|_{x=0}$,
which yields the identity $w=f_+'(0,\cdot)-m_-f_+(0,\cdot)$ in \er{T23-1}.
Substituting \er{T22-2} into $w=f_+'(0,\cdot)-m_-f_+(0,\cdot)$ we get \er{T23-1}.

Definitions of $s$ and $f_-$ give $s=\{f_+,\p_+\}=f_+{\p_+}'-{f_+}'|_{x=0}$,
which yields the identity $s=f_+(0,\cdot)m_+-f_+'(0,\cdot)$ in \er{T23-2}.
Substituting \er{T22-2} into the last identity we obtain \er{T23-2}.

These identities and analyticity of the functions
$f_+'(0,z), f_+(0,z), m_\pm(z)$
in the domain $\cZ$ imply that the functions
$\x(z), s(z), w(z), a(z)$ have analytic continuations
 from $z\in \R\sm \cup \ol g_n$ into $z\in\cZ$.
 The functions $a,b, s,w $ are analytic in $\cZ$ and are real on $i\R$. Then
the symmetry principle yields \er{symab}.

Asymptotics from Lemma \ref{T22}, \er{efs1}  and \er{Tm-1} and
identity \er{T23-1}, \er{T23-2} and  $\x=\vp(1,\cdot)w$
(see \er{iaw}) imply \er{T23-3}.

ii) We show $\r_-=2t$.  Due to \er{T23-3}, $s$ has exponential type $\r_-\le 2t$.
The decompositions $f_+=e^{ixz}(1+f)$ and $\p_+=e^{ixz} (1+\p)$
give $(1+f)(1+\p)=1+T, \ T=f+\p +\p f$ and
\[
\lb{esff}
s(z)=\int_0^tq(x)\p_+(x,z)f_+(x,z)dx=\int_0^tq(x)e^{i2xz}(1+T(x,z))dx,
\qqq  \qq z\in\cZ_\ve.
\]

Asymptotics \er{f2}, \er{efs1}, \er{T22-3} and $k(z)=z+O(1/z)$ as
$|z|\to \iy$ (see \er{pk1}) yield
\[
\lb{esff1}
\p(x,z)=O(1/z),\qq \qqq f(x,z)=e^{2(t-x)|\Im z|}O(1/z)\qq as \
 |z|\to \iy,\qq z\in \cZ_\ve.
\]
We need the following variant of the Paley-Wiener type Theorem from \cite{Fr}:

\no {\it Let $q\in \cQ_t^2$ and let each $G(x,z), x\in [0,t]$ be
analytic in $z\in\C_-$ and $G\in L^2((0,t)dx,\R dz)$.
Then $\int_0^te^{2izx}q(x)(1+G(x,z))dx$ has exponential type
at least $2t$ in $\C_-$.}

We can not apply this result to the function  $T(x,z), z\in\C_-$,
since $m_+(z)$ may have a singularity at $\m_n-i0\in \ol g_n+i0$ if
$g_n\ne \es$.  But we can use this result for the function
$T(x,z-i), z\in\C_-$, since \er{esff}, \er{esff1} imply $\sup _{x\in
[0,1]}|T(x,-i+\t)|=O(1/\t)$ as $\t\to \pm\iy$. Then the function
$s(z)$ has exponential type $2t$ in the half plane $\C_-$. The proof
for $\r_+=0$ is similar. \BBox

\section {Properties of the function $\x$}
\setcounter{equation}{0}

We start with the basic properties of the function $\x$.

\begin{lemma}
\lb{T32}
i) The function $\x=2ia(z)\sin k(z), z\in \cZ$ is analytic on
the Riemann surface $\cM$ and satisfies:
\begin{multline}
\lb{T32-1}
\x(z)=2i\sin k(z)(1+A(z))+J(z),\qq z\in \cM,\qqq \\
A(z)=\int_\R q(x)Y_2(x,z)dx,\qqq Y_2={1\/2}(\vp \wt\vt-\vt \wt\vp),\\
J(z)=-\int_\R q(x)Y_1(x,z)dx,\qqq Y_1=\vp(1,\cdot)\vt
\wt\vt-\vt'(1,\cdot)\vp \wt\vp+\b(\vp \wt\vt+\vt \wt\vp)
=\vp(1,\cdot,\cdot)+Y_{11},\\
Y_{11}=\vp(1,\cdot)\vt \wt\vt_*-\vt'(1,\cdot)\vp \wt\vp_*+
\b(\vp \wt\vt_*+\vt \wt\vp_*),
\ \
\wt\vt_*=\wt\vt-\vt, \wt\vp_*=\wt\vp-\vp,
\end{multline}
\begin{multline}
\lb{T32-2}
\x(z)=2(-1)^{n+1}(1+A(z))\sinh v(z)+J(z),\qqq z\in g_n^\pm\ne \es,\\
Y_{11}(z,x)=O(|g_n|/n^2)\ \ as \ \ z\in g_n, n\to \iy.
\end{multline}
where $v=\Im k$  and  $\pm v(z)>0$ for  $z\in g_n^\pm$.
The functions $J, A$ are entire and satisfy
\[
\lb{T32-30}
\x(z)=\ol \x(-\ol z), \qqq \qq \forall \ z\in \cZ,
\]
\[
\lb{T32-40}
J(z)=J(-z)=\ol J(\ol z)=\ol J(-\ol z), \qqq
 A(z)=A(-z)=\ol A(-\ol z)=\ol A(\ol z),\qq \forall \ z\in \C.
\]
\no ii)  The following estimates hold true
\begin{multline}
\lb{T32-3}
|J(z)|\le C_{p,q}\|q\|_t \rt(|\vp(1,z)|+
{|\vt'(1,z)|\/|z|^2}+{|\b(z)|\/|z|}\rt)e^{2t|\Im z|}
\le {C_0\/4|z|}e^{(2t+1)|\Im z|},
\\
|A(z)|\le {\|q\|_t^2C_{p,q}\/|z|^2}e^{2t|\Im z|},\qq where
\qq C_0=12\|q\|_te^{\|p\|_1+\|q\|_t+2\|p\|_t}, \qq
C_{p,q}=e^{2\|p\|_t+\|q\|_t},
\end{multline}
\[
\lb{T32-4}
|J(e_n)|\le {C_0\/|z_n|}\sinh h_n,\qqq if \qq e_n^-\ge 8e^{\|p\|_1}.
\]
\end{lemma}
\no {\bf Proof.} i)  Using \er{abws}, we rewrite the identity
\er{T23-1} in the form
\[
\lb{318}
\x(z)=2i\sin k(z)-\int_\R q(x)Y(x,z)dx,\qq
Y=\vp(1,\cdot)\p_-(x,\cdot)f_+(x,\cdot),
\]
for $z\in \C_+$. Let  $\vp_1=\vp(1,\cdot), \vt_1'=\vt'(1,\cdot)$. Using \er{318}, \er{T22-1} we rewrite $Y$ in the form
$$
Y=\vp_1(\vt+m_-\vp)(\wt\vt+m_+\wt\vp)=\vp_1\rt(\vt \wt\vt
+m_+\vt \wt\vp+m_-\vp \wt\vt-{\vt_1'\/\vp_1}\vp \wt\vp\rt)
=Y_1-i2Y_2\sin k.
$$
Substituting the last identity into $\int_\R q(x)Y(x,z)dx$ and using
\er{if} we obtain $Y_1=\vp(1,\cdot,\cdot)+Y_{11}$, which gives
\er{T32-1}.

\er{35} implies the identity in \er{T32-2}.
Substituting asymptotics from \er{efs}, Lemma \ref{T31} into $Y_{11}$
we obtain  $Y_{11}(z,x)=O(|g_n|/n^2)$ as $z\in g_n, n\to \iy.$

The function $\x$ is real on $i\R$, then the symmetry principle yields
\er{T32-30}.

The functions $A, J$ are entire and are real on $i\R, \R$. Then
 the symmetry principle yields \er{T32-40}.
Then $\x(z)$ is analytic  in $\cM$, since $\sin k(z)$ is analytic  in $\cM$

 ii) Using \er{efs}, \er{efs1}  and  \er{T32-1} and Lemma \ref{T31}, we obtain \er{T32-3}. Estimates \er{T32-3} and Lemma \ref{T31} give $|J(e_n)|\le {C_0\/2|z_n|} |g_n|$;
and the estimate $|g_n|\le 2h_n\le 2\sinh h_n$ (see \er{pk5}) yields \er{T32-4}.
\BBox

Define the functions $F, S$ by
$$
F(z)=|\x(z)|^2=4|\sin k(z)|^2|a(z)|^2>0, \qqq S(z)=|\vp(1,z)s(z)|^2,
\qq z\in \s_M=\R\sm \cup [e_n^-,e_n^+].
$$

\begin{lemma}
\lb{T33}
i)  The functions $F(z), S(z), z\in \s_M$
have analytic continuations into the whole complex plane $\C$ and satisfy
\[
\lb{T33-1}
F(z)=\x(z)\x(-z)=\x(z)\ol \x(\ol z),\qqq S(z)=\vp^2(1,z)s(z)s(-z),
\qqq z\in \cZ,
\]
\[
\lb{T33-11}
F=4(1-\D^2)(1+A)^2+J^2=4(1-\D^2)+S.
\]
Moreover, $F(z)>0$ and $S(z)\ge 0$ on each interval
$(e_{n-1}^+,e_n^-), n\ge 1$ and $F$ has even number of zeros on each
interval $[e_n^-,e_n^+], n\ge 1$.

\no ii) The function $F$ has only simple zeros at $e_n^\pm, g_n\ne \es$.
 Furthermore,

 if $e_n^-=e_n^+$ for some $n\ne 0$, then $e_n^-$ is
a double zero of $F$ and $e_n^-$ is not a state of $H$,

if $F(0)=0$, then $\z=0$ is
a double zero of $F$ and $\z=0$ is a virtual state of $H$.

\no iii) Let $\z\in \C_-\sm i\R$. The point $\z$ is a zero of $F$
iff $\z$ is a zero of $\x$ (with the same multiplicity).

\no iv) $\z\in i\R_-$ is a zero of $F$ iff $\z\in i\R_-$ or
$-\z\in i\R_+$ is a zero of $\x$.

\no v) Let $g_n\ne \es, n\ge 1$. The point $\z\in g_n$ is a zero of
$F$ iff $\z+i0\in g_n^+$ or $\z-i0\in g_n^-$ is a zero of $\x$ (with
the same multiplicity).

\end{lemma}

\no{\bf Proof.} i) Using  \er{T32-30} we deduce that $F=\x(z)\x(-z),
z\in \s_M$. Then by \er{T32-1}, $F$ satisfies
$F=(1-\D^2)(2-A)^2+J^2$ and then $F$ is entire. Using \er{iab},
\er{symab}, we obtain that $F, S$ satisfy \er{T33-1}, \er{T33-11}
and then $S$ is entire.

Recall that $F> 0$ and $S\ge0$ inside each $(e_{n-1}^+,e_n^-)$,
since  $|\sin k(z)|>0$ and  $|a(z)|\ge1$ on $\R\sm \cup_{n\ne 0} \ol
g_n$, see \er{iab}. Due to $F(e_n^\pm)\ge 0$, we get that $F$ has
even number of zeros on each interval $[e_n^-,e_n^+]$.

ii) Consider the case $\z=e_n^+, g_n\ne \es$, the proof for
$\z=e_n^-$ is similar. We have $F=F_0+S$, where $F_0=4(1-\D^2)$.
Thus if $F(\z)=0$, then we get $S(\z)=0$, since $\D^2(e_n^+)=1$.
Moreover, $F_0'(\z)=-2\D(\z)\D'(\z)>0$ and $S'(\z)\ge 0$, which
gives that $\z=e_n^+$ is a simple zero of $F$.

Let $g_n=\es$. Then the functions $\vp_1, \vt_1'$ have zero at
$e_n=e_n^\pm$ and then identity \er{LD0} gives $\b(e_n)=0$. Then
identity \er{T32-1} gives $J(e_n)=0$ and identity \er{T33-11}
implies that $e_n$ is a zero of $S$ with the multiplicity $\ge 2$.
Differentiating  $F=4(1-\D^2)+S$  we obtain $\ddot F(e_n)=-8\ddot
\D\D+\ddot S|_{e_n}\ge -8\ddot \D\D|_{e_n}>0$, since $\ddot
S(e_n)\ge 0$. Then $e_n$ is a second order zero of $F$.

The proof for the case $F(0)=0$ is similar.

iii) The function $\x$ has not zeros in $\C_+\sm i\R$, see \cite{F1}.
This and the identity $F=\x(z)\x(-z)$ yields v).

iv) The definition $F(z)=\x(z)\x(-z), z\in \cZ$ gives iv).

v)  The statement v) follows from iii).
\BBox

Due to Lemma \ref{T33} we study the entire function $F$ instead of
the function $\x$ on the Riemann surface
$\cM$. Now we describe the forbidden domain for the resonances.

\begin{lemma}
\lb{T44}
$F$ has not zeros in $\cD_0\sm \cup [e_n^-,e_n^+]$, where
\[
\lb{dFD}
\cD_0=\rt\{z\in \C: |z|>\max \{180e^{2\|p\|_1},C_0e^{2t|\Im z|}\}\rt\},\qq
C_0=12\|q\|_te^{\|p\|_1+\|q\|_t+2\|p\|_t}.
\]
Moreover, if $[e_n^-,e_n^+]\ss \cD_0$, then  $F$ has exactly two
zeros $z_n^\pm\in [e_n^-,e_n^+]$ such that:

if $e_n^-<e_n^+$, then $z_n^-, z_n^+$ are the simple zeros such that
$e_n^-\le z_n^-<e_n<z_n^+\le e_n^+$,

if $e_n^-=e_n^+$, then $z_n^-=z_n^+$ is a zero of order two,

There are no other zeros of $F$ in $\cD_0$.

\end{lemma}

\no{\bf Proof.} Using Lemma \ref{T32}, \er{efs} and \er{efs1}  we obtain
$$
|A(z)|\le {\|q\|_t^2\/|z|_1^2}X_2,\qqq |J(z)|\le {3\|q\|_t\/|z|_1}
X_2X\qqq |\D(z)-\cos z|\le {X\/|z|_1},\qq z\in \C,
$$
where $ X_2=e^{|\Im z|2t+\|q\|_t+2\|p\|_t},\qq X=e^{|\Im
z|+\|p\|_1}$ and $|z|_1=\max\{1,|z|\}$. Substituting these estimates
into the identity
$$
F-4\sin^2 z=4(\cos^2z-\D^2)+J^2+4(1-\D^2)A(A+2),
$$
which follows from \er{T33-1},  we obtain
$$
|F(z)-4\sin^2 z|\le 9X^2C_1,\qqq
C_1={1\/|z|}+{\|q\|_t^2X_2^2\/|z|^2}+{\|q\|_t^2X_2\/|z|^2}
\rt(2+{\|q\|_t^2X_2\/|z|^2}\rt),\qq |z|\ge 1.
$$
Using the simple estimate $e^{|{\Im}z|}\le 4|\sin z|$
for all $|z-\pi n|\ge {\pi \/4}, n\in \Z$, (see p. 27 [PT]), we obtain
$$
9X^2=9e^{2|\Im z|+2\|p\|_1}\le  |4\sin^2 z|r_0^2
 \qq all \ |z-\pi n|\ge{\pi\/4},\qq  n\in \Z, \qqq r_0=6e^{\|p\|_1},
$$
which yields
$$
|F(z)-4\sin^2 z|\le 4|\sin^2 z|r_0^2C_1< 4|\sin^2 z|,\qq \forall
\qq z\in \cD_0\sm \bigcup \dD_{\pi\/4}(\pi n),
$$
since for $z\in \cD_0$ the following estimates hold true
$$
{r_0^2\/|z|}<{1\/5},\qqq\qqq {r_0\|q\|_tX_2\/|z|}<{1\/2},
$$
$$
r_0^2C_1\le
{r_0^2\/|z|}+{r_0^2\|q\|_t^2X_2^2\/|z|^2}+
{r_0^2\|q\|_t^2X_2\/|z|^2}\rt(2+{\|q\|_t^2X_2\/|z|^2}\rt)<{1\/5}
+{1\/4}+{1\/2}+{1\/(24)^2}<{19\/20}.
$$
Thus  by Rouche's theorem, $F$ has as many roots, counted with
multiplicities, as $\sin^2 z$ in each disk $\dD_{\pi\/4}(\pi n)\ss\cD_0$.
Since $\sin z$ has only the roots $\pi n, n\ge 1$, then $F$
has two zeros in each disk $\dD_{\pi\/4}(\pi n)\ss\cD_0$ and $F$ has
not zeros in $\cD_0\sm \cup \dD_{\pi\/4}(\pi n)$.

Let $e_n^-<e_n^+$.  Estimate \er{T32-3} gives $|A(e_n)|\le {1\/2}$
and \er{T32-4} imply $|J(e_n)|\le \sinh h_n$. The substitution of
these estimates into \er{T33-1} yields
 \[
\lb{T33-2}
F(e_n)\le - \rt(1-{C_0^2\/|e_n|^2}  \rt)\sinh^2 h_n<0,
\]
and $e_n,h_n$ are defined by \er{dzn}. The function $F(e_n^\pm)\ge
0$ and due to \er{T33-2}, we deduce that $F$  has exactly  two zeros
$z_n^\pm $ on the segment $[e_n^-,e_n^+]$ such that $e_n^+\le
z_n^-<e_n<z_n^+\le e_n^+$.

If $e_n^-=e_n^+$, then by Lemma \ref{T33} ii),   $z_n^-=z_n^+$ is a
second order zero of $F$. \BBox

\bigskip

We discuss the relationship of states from the Definition S and poles
of the resolvent $(H-z^2)^{-1}$.
The kernel of the resolvent $R=(H-z^2)^{-1},  z\in \C_+,$ has the form
$$
R(x,x',z)={f_- (x,z )f_+(x',z)\/-w(z)}={R_1(x,x',z)\/-\x(z)},\ \ \ x<x',\ \
R_1=\vp(1,z)f_- (x,z )f_+(x',z),
$$
and $R(x,x',z )=R(x',x,z ),\ x>x'$. Identity \er{T22-1} yields
$f_\pm=\wt\vt+m_\pm \wt\vp,\qq \wt\vt=\wt\vt(x,z),\qq  \wt\vp=\wt\vp(x,z)$. Let $ \wt \vt_*=\wt\vt(x',z),\qq \wt \vp_*=\wt\vp(x',z)$. Then using \er{LD0} $m_\pm={\b\pm i\sin k\/\vp(1,\cdot)}$,  we obtain
$$
R_1(x,x',z)=\vp(1,\cdot)\wt \vt \wt \vt_*+(\b-i\sin k)\wt \vp \wt \vt_*+
(\b+i\sin k)\wt\vt \wt \vp_*-\vt'(1,\cdot)\wt \vp_*\wt \vp.
$$
Then for fixed $x,x'\in \R$  the function $R_1(x,x',z)$ is analytic
on $\cM$ and $R_1$ is locally bounded on $\R^2\ts \cM$. The zeros of
$\x$ create the singularities of the kernel $R(x,x',z)$. Thus if
$\x(z)=\vp(1,z)w(z)=0$ at some $z\in\cM$, then  $R(x,x',z)$ has
singularity at $z$. The poles of $R(x,x',z)$ define the bound states
and resonances. The zeros of $\x$ define the bound states and
resonances, since the function $R_1$ is locally bounded.

Consider the unperturbed case $q=0$. Recall that $\sin k(z)$ is
analytic  in $\cM$ and $\sin k(z)=0$ for some $z\in \cM$ iff
$z=e_n^-$ or $z=e_n^+$  for some $n\ge 0$.  The function $\x$ is
analytic on $\cM$ and has branch points $e_n^\pm, g_n\ne \es$. Then
$R_0=(H_0-z^2)^{-1}$ has the form
$$
R_0(x,x',z)={R_{10}(x,x',z)\/-\x_0(z)},\qq \x_0=\vp(1,z)w_0(z)=2i\sin k(z),
$$$$
R_{10}=\vp(1,\cdot)\p_-(x,\cdot)\p_+(x',\cdot)=
\vp(1,\cdot)\vt\vt_*+(\b-i\sin k)\vp\vt_*+
(\b+i\sin k)\vt\vp_*-\vt'(1,\cdot)\vp\vp_*,
$$
where $\vp=\vp(x,\cdot),..$ and $\vp_*=\vp(x',\cdot),..$ Thus
$R_0(x,x',z)$ has singularity at some $z\in\cM$ iff $\sin k(z)=0$,
i.e., $k(z)=\pi n $ and then $z=e_n^\pm$.

Remark that if  $\z=e_n^-=e_n^+$ (i.e.. the gap $g_n=\es$), then
$\x$ is analytic at  $\z$, and  $\x(\z)=0$, but such point $\z$ is
not the state. The function $a(z)$ is analytic at $\z$ and \er{iab}
yields  $|a(\z)|\ge 1$

\section {Proof of Theorems 1.1-1.4}
\setcounter{equation}{0}

{\bf Proof of Theorem \ref{T1}}.
i) By Lemma \ref{T32}, $\x$ is analytic on $\cM$ and the function $J$
is entire.

ii) Recall that \er{T33-1} gives
$$
F(z)=\x(\z)\x(-\z)=\x(\z)\ol \x(\ol \z), \qqq \z=z+i0\in g_n+i0.
$$
Then the zeros $\z\in g_n+i0$ of $\x(\z)$ give the bound states and
the zeros $\z\in g_n+i0$ of $\x(\ol\z)$ give the antibound states.
Their global number on $g_n^c$ plus the possible virtual states at
$e_n^\pm$ is even $\ge 0$, see Lemma \ref{T33}, i).

The similar arguments and Lemma \ref{T44} yield iii).

Moreover, if an open gap $g_n=(e_n^-,e_n^+)\ss \cD_0$ ( $\cD_0$ is defined by \er{T44}), then  there exist exactly two simple zeros $z_n^\pm\in [e_n^-,e_n^+]$ such that
\[
\lb{5e}
e_n^-\le z_n^-<e_n<z_n^+\le e_n^+.
\]
The asymptotics \er{aen} yields $\d_n^\pm \le {2\/3}|g_n|$ as $n\to \iy$.
Note that if $g_n=\es$, then $F$ has a double zero $e_n^\pm=z_n^\pm=e_n$.
There are no other zeros of $F$ in $\cD_0$.

iv) Due to \er{5e} we have  $z_n^\pm=e_n^\pm\mp \d_n^\pm,
\d_n^\pm\ge 0$. Let $\z=z_n^\pm,  \d=\d_n^\pm$. Then the equation
$0=\x(\z)=(-1)^{n+1}2(1+A(\z))\sinh v(\z)+J(\z),\  \ \z\in \ol
g_n^\pm\ne \es$ and Lemma \ref{T31}, \er{T32-3} imply
$$
\sinh |v(\z)|=O(J(\z)) =O\rt(
|\vp(1,\z)|+{|\vt'(1,\z)|\/|\z|^2}+{|\b(\z)|\/\z}\rt) =\ve_n
O(|g_n|), \qq \ve_n={1\/2\pi n}
$$
as $n\to \iy$.  Moreover, using the estimate
$|(z-e_n^-)(z-e_n^+)|^{1\/2}\le |v(z)|$ for each $z\in g_n$ (see
\er{prq}) we obtain $|\d(|g_n|-\d)|^{1\/2}\le |v(\z)|=\ve_n
O(|g_n|)$, which yields $\d=\ve_n^2 O(|g_n|)$. Thus the points
$z_n^\pm$ are close to $e_n^\pm$ and satisfy:
\[
\lb{asd} |v(z_n^\pm)|=\ve_n O(|g_n|),\qqq\qqq
\d_n^\pm=z_n^\pm-e_n^\pm=\ve_n^2 O(|g_n|).
\]

Consider the case $\d=\d_n^-=\ve_n^2O(|g_n|)$, the proof for
$\d=\d_n^+$ is similar. Using \er{T32-2} we obtain
$$
J=J_{0}+O(\ve_n^2 |g_n|),\qqq J_{0}(z)=-\int_\R\vp(1,z,x)q(x)dx.
$$
Asymptotics \er{dex} gives
\[
\lb{J0} J_{0}(\z)=(-1)^n|g_n|\ve_n (\wh q_0-c_n\wh q_{cn}+s_n\wh
q_{sn}+ O(\ve_n)),
\]
where $c_n=\cos \f_n(0), s_n=\sin \f_n(0)$, and \er{T32-2} yields
\[
\lb{J-n} (-1)^{n}J(\z)=\ve_n |g_n|I_n^-, \qqq I_n^-=\wh q_0-c_n\wh
q_{cn}+s_n\wh q_{sn}+ O(\ve_n).
\]
Using \er{T32-3} we obtain
\[
\lb{as2} \sinh v(\z)={(-1)^{n}J(\z)\/2+2A(\z)}={\ve_n
|g_n|I_n^-\/2+O(\ve_n^2)},\qq \sign v(\z)=\sign (-1)^{n}J(\z)=\sign
I_n^-.
\]
Note  that if $v(\z)>0$ , then $\z\in g_n^+$ is a bound state,

 if $v(\z)<0$,  then $\z\in g_n^-$ is a resonance,

if $v(\z)=0$, then  $\z=e_n^-$ or $\z=e_n^+$ is a virtual state.

Then  \er{asd} gives $\sinh v(\z)=v(\z)(1+O(|g_n|^2\ve_n^2))$  and
using asymptotics \er{pav}, we obtain
$$
 v(\z)=\sqrt{\d(|g_n|-\d)}(1+O(\ve_n^2))=\sqrt{\d|g_n|}(1+O(\ve_n^2)).
$$
This and \er{as2} yield $\d_n^-={|g_n|\ve_n^2\/4}(I_n^-)^2$ and
\er{ape} gives $|\g_n|=(2\pi n)|g_n|(1+O(\ve_n^2))$, which yields
$\d_n^-={2|\g_n|\/(4\pi n)^3}(I_n^-)^2(1+O(\ve_n^2))$. This and
\er{J-n} give \er{T2-1ad}.

If $\wh q_0>0$, then $I_n^->0$ and  above arguments yield that
$z_n^-$ is a bound state and $z_n^+$ is an antibound state.
Conversely, if $\wh q_0<0$, then $I_n^-<0$ and we deduce that
$z_n^-$ is an antibound state and $z_n^+$ is a bound state. This
yields \er{q0}.\BBox

Note that due to \er{T2-1ad}, the high energy real states of $H$ and
$H_0$ are very close.
 This gives
 \[
\lb{T2-3}
\#(H,r,\cup_{n\ge 1} g_n^c)=\#(H_0,r,\cup_{n\ge 1} g_n^c)
+2N_* \qq as \qq r\to \iy, \qq r\notin \cup_{n\ge 1}\ol  g_n.
\]
for some $N_*\in \Z$.

\no {\bf Proof of Corolarry \ref{Tc}}.
We have $-c_n\wh q_{cn}+s_n\wh q_{sn}=-|\wh q_n|\cos (\f_n+\t_n)$, where
$c_n, s_n$ given by \er{cs}. Then
Theorem \ref{T1} iv) yields the Statement i) and ii).

If $p\in L_{even}^2(0,1)$, then the coefficient $s_n=0$ for all $n\ge 1$
(see remark before
 Corolarry \ref{Tc}). Thus Statement i) yields Statement iii).
\BBox

\no {\bf Proof of Theorem \ref{T2}}.  An entire function $f(z)$ is
said to be of $exponential$ $ type$  if there is a constant $\a$
such that  $|f(z)|\leq $ const. $e^{\a |z|}$ everywhere. The
function $f$ is said to belong to the Cartwright class $Cart_\r,$ if
$f$ is entire, of exponential type, and the following conditions
hold true:
$$
\int _{\R}{\log (1+|f(x)|)dx\/ 1+x^2}<\iy  ,\ \
\r_\pm(f)=\r,\ \ \ {\rm where}\ \ \
\r_{\pm}(f)\ev \lim \sup_{y\to \iy} {\log |f(\pm iy)|\/y}.
$$

Let $\cN (r,f)$ be the total number of  zeros of $f$  with modulus $\le r$,
each zero being counted according to its multiplicity. We recall the
well known result (see [Koo]).

\no    {\bf Levinson Theorem.} {\it  Let the entire function $f\in
Cart_\r,\r>0$. Then $ \cN(r,f)={2r\/ \pi }(\r+o(1))$ as $r\to \iy$,
and for each $\d >0$ the number of zeros of $f$ with modulus $\le r$
lying outside both of the two sectors $|\arg z | , |\arg z -\pi
|<\d$ is $o(r)$ for $r\to \iy$.}

Consider the functions $F, S$.

By Lemma Lemma \ref{T33}, the functions $F, S$ are entire and by
\er{T32-3}, the function $F\in L^\iy(\R)$ and then the function
$S\in L^\iy(\R)$ . Using  Lemma \ref{T23} ii) and \er{efs}, we
deduce the function $S$ has the exponential type $\r_\pm(S)=2+2t$ in
the half plane $\C_\pm$. Thus the function $S$ belongs to
$Cart_{2+2t}$ and the identity \er{T33-11} gives that $F\in
Cart_{2+2t}$ and the Levinson Theorem implies
\[
\lb{56}
\cN(r,F)=2r{2+2t+o(1)\/\pi}\qq as \qq r\to \iy.
\]

Let $\pm \z_n>0, n\ge 1$  be all real zeros $\ne 0$  of $F$ and let
the zero $\z_0=0$ have the multiplicity $n_0\le 2$. Define the
function $F_1=z^{n_0}\lim_{r\to \iy}\prod_{|\z_n|\le
r}(1-{z\/\z_n})$. Recall that by Lemma \ref{T33}, $F(z)>0$ on the
set $\R\sm \cup_{n\ne 0} g_n$ and by Lemma \ref{T44}, the function
$F$ has exactly two zeros  on each set $[e_n^-,e_n^+]$ for  $n$
large enough. Then
\[
\lb{57}
\cN(r,F_1)=2r{2+o(1)\/\pi} \ \ \as \ r\to\iy.
\]
Combining \er{56} and \er{57}, we obtain
$$
\cN(r,F/F_1)=\cN(r,F)-\cN(r,F_1)=2r{2+2t+o(1)\/\pi}-2r{2+o(1)\/\pi}
=4r{t+o(1)\/\pi}.
$$
Denote by $\cN_+(r,f)$ (or $\cN_-(r,f)$) the number of zeros of $f$
with imaginary part $>0$ (or $<0$) having modulus $\le r$, each zero
being counted according to its multiplicity.
 The function $F$ is real on the real line,  then
 $$
 \cN_+(r,F)=\cN_-(r,F)={1\/2}\cN(r,F/F_1)=2r{t+o(1)\/\pi}.
 $$
Then Lemma \ref{T33} gives the identities  $\cN_-(r,F)=\cN_+(r,F)=\cN_-(r,\x)+N_*$ for some integer $N_*\ge 0$,
which yields $\cN_-(r,\x)=2r{t+o(r)\/\pi}$ as $r\to \iy$ and \er{T2-2}.
\BBox

\no {\bf Proof of Theorem \ref{T3}}. i) Let the operator $H_0$ have
infinitely many gaps $\g_n\ne \es$ for some $p\in L^2(0,1)$ and let
$\vk=(\vk_n)_{1}^\iy$ be any sequence, where $\vk_n\in \{0,2\}$. For
this $p$ there exist a unique sequence of angles $\f_n\in [0,2\pi),
n\ge 1$, defined by  \er{cs}.

We take a real potential $q\in L^2(0,1), \supp q\ss (0,1)$ given by
\[
\lb{defq}
q(x)=\sum_{n\ge 1}{1\/|n|^\a}(e^{i\t_n+i2\pi nx}+e^{-i\t_n-i2\pi nx}),
\qq x\in (0,1), \qq {1\/2}<\a<1.
\]
Let $\ve\in (0,1)$. We take $\t_n\in [0,2\pi )$ such that

if   $\vk_n=2$, then we choose  $\t_n$ such that  $\cos (\f_n+\t_n)<-\ve$,

if $\vk_n=0$, then we choose  $\t_n$ such that  $\cos (\f_n+\t_n)>\ve $.

Then due to Corolarry \ref{Tc} i), the operator $H$ has
$\vk_n=1-\sign \cos (\f_n+\t_n)$ bound states in the physical gap
$g_n^+\ne \es$ and $2-\vk_n$ resonances inside the  nonphysical gap
$g_n^-\ne \es$  for $n$ large enough.

\medskip

ii) Let  $q\in \cQ_t,t>0$ satisfy  $\wh q_0=0$ and let $\wh
q_{n}=|\wh q_{n}|e^{i\t_n}$, where $|\wh q_{n}|>n^{-\a}$ and
$\t_n\in [0,2\pi)$  for all $n$ large enough  and some $\a\in
(0,1)$. Let $\vk=(\vk_n)_{1}^\iy$ be any sequence, where $\vk_n\in
\{0,2\}$.

Let $\d=(\d_n)_1^\iy\in \ell^2$ be a sequence of nonnegative numbers
$\d_n\ge 0, n\ge 1$ and infinitely many $\d_n>0$ and let
 $(\wt\f_n)_1^\iy$ be a sequence of angles $\wt\f_n\in [0,2\pi), n\ge 1$.
Let $\ve\in (0,1)$. We take $\wt\f_n\in [0,2\pi )$ such that

if   $\vk_n=2$, then we choose  $\wt\f_n$ such that  $\cos (\wt\f_n+\t_n)<-\ve$,

if $\vk_n=0$, then we choose  $\wt\f_n$ such that  $\cos (\wt\f_n+\t_n)>\ve $.

Recall the result from \cite{K5}:

{\it The mapping $\P: \cH\to \ell^2\os \ell^2$ given
by $\P=((\P_{cn})_1^\iy,(\P_{sn})_1^\iy)$ where
\[
\lb{ip1}
 \P_{cn}={E_n^-+E_n^+\/2}-\m_n^2,\qq
\P_{sn}= \rt|{|\g_n|^2\/4}-\P_{cn}^2\rt|^{1\/2}\sign (|\vp'(1,\m_n)|-1),
\]
is a real analytic isomorphism between real Hilbert spaces
$\cH=\{p\in L^2(0,1): \int_0^1p(x)dx=0\}$ and $\ell^2\os \ell^2$. }
Note that \er{ip1} and \er{cs} give $\P_{cn}={|\g_n|\/2}c_n,
\P_{sn}={|\g_n|\/2}s_n$, since $\vp(\cdot,\m_n)=y_n(\cdot)$.

Then for any sequence $\d=(\d_n)_1^\iy\in \ell^2$ and any sequence
of angles $(\wt\f_n)_1^\iy$ there exists a unique potential $p\in
\cH\ss L^2(0,1)$ such that each gap length  $|\g_n|=\d_n$  and the
corresponding angle $\f_n=\wt\f_n$ for all $n\ge 1$. We consider the
operator $H=-{d^2\/dx^2}+p+C+q$, where  the constant $C$ is such
that $E_0^+=0$.

Then due to Corolarry \ref{Tc} i), the operator $H$ has
$\vk_n=1-\sign \cos (\f_n+\t_n)$ bound states in the physical gap
$g_n^+\ne \es$ and $2-\vk_n$ resonances inside the  nonphysical gap
$g_n^-\ne \es$  for $n$ large enough. \BBox

\

\no {\bf Acknowledgments.} \small The various parts of this paper
were written at  ESI, Vienna, and Mathematical Institute of the
Tsukuba Univ., Japan and Ecole Polytechnique, France. The author is
grateful to the Institutes for the hospitality.
This work was supported by the Ministry of education and science of the
Russian Federation, state contract 14.740.11.0581.


\end{document}